# Global Search of Optimal Spacecraft Trajectories using Amortization and Deep Generative Models

Ryne Beeson[*], Anjian Li[†], and Amlan Sinha[‡]
*Princeton University, New Jersey, 08540*

**Preliminary spacecraft trajectory optimization is a parameter dependent global search problem that aims to provide a set of solutions that are of high quality and diverse. In the case of numerical solution, it is dependent on the original optimal control problem, the choice of a control transcription, and the behavior of a gradient based numerical solver. In this paper we formulate the parameterized global search problem as the task of sampling a conditional probability distribution with support on the neighborhoods of local basins of attraction to the high quality solutions. The conditional distribution is learned and represented using deep generative models that allow for prediction of how the local basins change as parameters vary. The approach is benchmarked on a low thrust spacecraft trajectory optimization problem in the circular restricted three-body problem, showing significant speed-up over a simple multi-start method and vanilla machine learning approaches. The paper also provides an in-depth analysis of the multi-modal funnel structure of a low-thrust spacecraft trajectory optimization problem.**

## Nomenclature

| | | |
|---|---|---|
| $O$ | = | optimal control problem |
| $J$ | = | objective function |
| $\phi$ | = | terminal cost |
| $\mathcal{L}$ | = | running cost |
| $\mathcal{L}_{\text{CVAE}}$ | = | CVAE loss function |
| $\mathcal{U}$ | = | function space for continuous-time control |
| $\xi$ | = | state solution |
| $f$ | = | uncontrolled state vector field |
| $g$ | = | controlled state vector field |
| $t_0$ | = | initial time |

[*]Assistant Professor, Department of Mechanical and Aerospace Engineering, Member AIAA
[†]Ph.D. Candidate, Department of Electrical and Computer Engineering
[‡]Ph.D. Candidate, Department of Mechanical and Aerospace Engineering, Member AIAA

| | | |
|---|---|---|
| $t_f$ | = | final time |
| $t$ | = | time |
| $\psi$ | = | path constraints |
| $\mathcal{K}$ | = | path constraint index set |
| $\Xi$ | = | optimal control boundary conditions |
| $\mathbb{H}$ | = | set of control transcriptions |
| $x$ | = | primal variables |
| $\mathcal{X}$ | = | primal space |
| $c$ | = | nonlinear programming constraints |
| $\mathcal{E}$ | = | equality index set |
| $\mathcal{I}$ | = | inequality index set |
| $\mathcal{P}$ | = | nonlinear program |
| $\pi$ | = | numerical parameter optimization solver |
| $\Lambda$ | = | dual space |
| $\lambda$ | = | dual variables |
| $\mathcal{U}^h$ | = | primal-dual space |
| $\mathcal{A}$ | = | collection of local optimal primal-dual solutions |
| $\hat{d}$ | = | diversity "metric" |
| $\eta$ | = | diversity threshold parameter |
| $p = p(\cdot\|\alpha)$ | = | conditional probability distribution |
| $\mathcal{N}(\cdot;\mu,\Sigma)$ | = | normal distribution with mean $\mu$ and covariance $\Sigma$ |
| $\mathcal{G}$ | = | Gaussian mixture model |
| $\mu_L$ | = | Lebesgue measure |
| $U$ | = | uniform distribution |
| $\delta$ | = | Dirac distribution |
| $^k\mathcal{N}$ | = | local neighborhood of $k$ iterations |
| $\tau$ | = | time variables |
| $m$ | = | mass variable |
| $u$ | = | thrust control vector |
| $z$ | = | latent variable |
| $q_\phi$ | = | posterior distribution (encoder) |
| $p_\theta$ | = | prior distribution |

| $p_\psi$ | = | likelihood (decoder) |

Superscripts

| $k$ | = | numerical solver iteration index |
| $h$ | = | parameter index for control transcription |
| $*$ | = | numerical solutions that are extrema |

Subscripts

| $\alpha$ | = | parameter index for optimal control problem |
| $h$ | = | parameter index for control transcription |
| $\gamma$ | = | parameter index for parameter optimization solver |
| $\beta$ | = | threshold for objective value |
| $i, j, k$ | = | indices |
| $\psi, \phi, \theta$ | = | hyperparameters for neural networks |

## I. Introduction

The preliminary spacecraft trajectory design phase can be posed as a parameterized global search problem for optimal spacecraft trajectories. At each stage of the preliminary design, the mission objectives, requirements, and constraints may change, resulting in variations of the global search problem parameters. Parameters may also change to represent increased modeling fidelity. The aim at any stage of the preliminary design is to solve for a large set of high quality spacecraft trajectories with diverse, or similarly qualitatively different, features. High quality is naturally defined by the value of a solution's objective value relative to the best known. Examples of qualitatively different features may include trajectories that have a different number of revolutions around a central body, a different number or sequence of gravity assist flybys, solutions that avoid radiation belts or other hazards, or solutions that depart the original or target orbital planes. The benefit of having different qualitative solutions is that it allows mission designers to trade different priorities in their design and reflects the fact that not all relevant objectives and constraints can be incorporated into the optimal spacecraft trajectory problem so early or readily in the design phase (i.e., without prior knowledge of what is relevant and when designing at a quick cadence).

In the simplest of cases, a mission designer's past experience may be sufficient to guide them in finding a high quality set of solutions. More often than not, experience and intuition fall short when designing complex trajectories. This is for instance typically true when designing for low-thrust (LT) propulsion systems, especially in multiple gravitational body regimes.

Due to the fact that closed form solutions rarely exist for practical continuous-time optimal control problems, a control transcription is used to convert either the indirect or direct formulation of the optimal control problem into a finite dimensional parameter optimization problem (see Sec. II.B for further discussion and references on indirect and direct transcriptions). The parameter optimization problem is then solved using an appropriate gradient based numerical solver (NS). The parameter optimization problem can be of high dimension (e.g., hundreds to tens of thousands of parameters), especially in the direct optimal control approach, and requires an initial guess for the NS. Performing an exhaustive grid search over this high dimensional problem is numerically intractable. In the case of an indirect optimal control approach, the dimensionality of the parameter optimization problem can typically be kept low, but the costate variables lack physical meaning and can have feasible domains with poor boundary regularity. The result of these difficulties on the parameterized global search problem is again a computational intractability if a naive exhaustive grid search is performed.

The global search process can be thought of as having three algorithm levels that work together to generate a large collection of useful solutions. A depiction of a global search algorithm for a fixed problem is given in Fig. 1. The first-level (highest-level or simply Level-1 in Fig. 1) attempts to provide initial guesses within the neighborhood of locally optimal solutions. In particular, the initial guesses should be within the basins of attraction to local optimal solutions (examples of basins are shown as gray columns in the second funnel region of Fig. 1). One can think of these initial guesses as coming from a probability distribution defined on the parameter optimization state space. The action of sampling the distribution (blue curve) is shown by the pink vertical dashed arrow in Fig. 1. When a global optimization problem has a monotonically decreasing set of objective values for local minima at some relatively large scale, that subset of the state space is said to have a "funnel" structure; a terminology adopted from chemists and physicists in the molecular conformation research community (see for instance Wales and Doye [1], Leary [2], or Wales [3]). It is not only the number of high quality local minima (e.g., local minima below the $\beta$-threshold shown in Fig. 1) that makes solution of global optimization problems difficult, but also the relative geometry of the local minima, and the topology of their basins of attraction and funnel structures. The problem considered in this paper is one in which this topology of solutions is parameterized (i.e., the parameterized global search problem) and therefore we aim to solve a family of these problems.

The third-level (lowest-level or Level-3 in Fig. 1) performs a local gradient based optimization given an initial guess to find a local optima. The initial guess may come from the first-level algorithm, or it may come from the second-level (intermediate-level and also labeled as Level-2 in Fig. 1) algorithm, which provides small perturbations of states to encourage traversal of the funnel structure towards the minima of the funnel.

Research in the astrodynamics community has largely focused on methods for solving a fixed parameterization of the global search problem and has made use of methods from the global optimization community such as the simple multi-start algorithm, heuristic methods that include genetic algorithms, particle swarms, differential evolution, and simulated annealing, as well as basin hopping variants [4]. Since our goal in this paper is learning the structure and topology of the local optimal solutions, their basins and neighborhoods, and the funnel structure as the global optimization problems change through an index set (i.e., parameterization), we restrict our literature review to work that does not simply apply global optimization techniques, but also attempts some hypothesis or observation of the structure and topology of these sets.

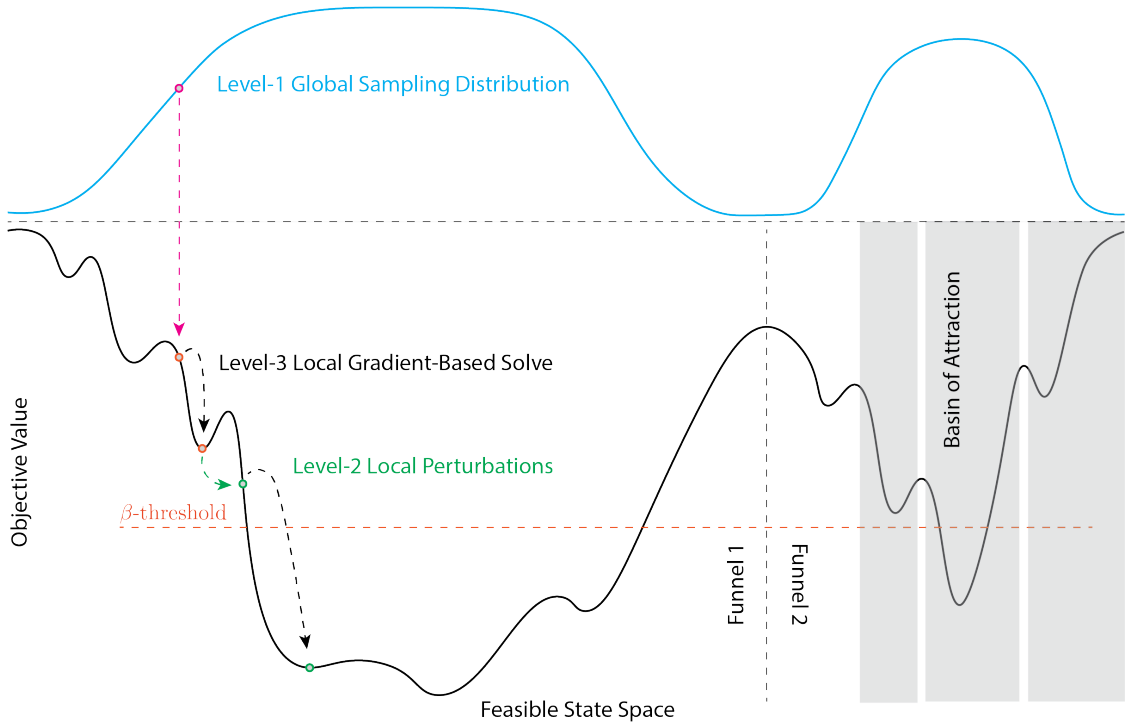

**Fig. 1  The three generic levels of a global search algorithm over a multiple funnel, multi-modal feasible state space.**

It was Hartmann et al. [5], building on the work of Rauwolf and Coverstone-Carroll [6], that first started exploring the global solution landscape of LT trajectories by using stochastic search methods. They produced visualizations of multiple objective Pareto fronts for an Earth-Mars direct transfer problem. In their work, they provided initial guesses to the indirect Solar Electric Propulsion Trajectory Optimization Program (SEPTOP) solver using a non-dominated genetic algorithm. The paper by Russell [7] demonstrated the local optima solution landscape for constant specific impulse (CSI) LT solutions solved with the (indirect) primer-vector theory [8] in planar circular restricted three-body problems (PCR3BP) of the Jupiter-Europa and Earth-Moon systems by also creating Pareto fronts; in Russell's paper, this was for $\Delta v$ (equivalently, delivered mass) versus time-of-flight. They made use of an exhaustive grid search for supplying initial guesses. Oshima et al. [9] provided an extended analysis of the Earth-Moon PCR3BP given in Russell. They analyzed the behavior of solutions with respect to a Tisserand-Poincaré graph. Although some idea of the global solution landscape can be inferred from Pareto fronts, it does not provide the topology of the solutions which belong to the control space; the visualization of Pareto fronts are in the function-valued space (e.g., $\mathbb{R}^2$ in Russell [7] and Oshima et al. [9], and $\mathbb{R}^3$ in Hartmann et al. [5]).

For the decade following the work of Hartmann et al. [5], efforts by the community to understand the topology of solutions in the control space was largely constrained to patched-conic approximations for multiple gravity assist (MGA) problems, sometimes involving deep space maneuvers (DSM). The technical reports of Di Lizia and Radice [10] and Myatt et al. [11] provided the first in-depth analysis of solutions in the control space for a number of simple transfers. These mostly included impulsive, some of MGA-DSM type, but a few planet-to-planet examples of LT were included. Constraint sets were always of the simple box-constrained type.

Making use of the MGA-DSM model, Vasile and De Pascale [12] showed that even the inclusion of a single DSM in a direct planet-to-planet transfer quickly increases the complexity of the feasible search domain, making it highly nonconvex and multimodal. They also show that it greatly increases the number of high quality local minima, as well as significantly alters the distribution of those minima in the control space. Vasile et al. [13] extended the study of MGA-DSM problems by applying an approach developed by Reeves and Yamada [14], which partitions the range space of the objective function into levels (i.e., a collection of $\beta$-thresholds as indicated in Fig. 1) and then calculates the intralevel and translevel distances between local optima. Visualizing these distances, one can then attempt to infer as to whether funnel structures are present and choose appropriate tuning parameters for global optimization algorithms. Vasile et al. [13] use the aforementioned analysis to develop a new variant of differential evolution with behavior closer to that of monotonic basin hopping (MBH); allowing it to more easily navigate to the bottom of funnels, but without the additional regularity requirements of MBH. Addis et al. [15] also studied MGA-DSM problems put forth by the European Space Agency Advanced Concepts Team (the Global Trajectory Optimization Problem database [16]), applying various global optimization algorithms and conjecturing that these problems display "a funnel structure similar, in some sense, to that of molecular optimization problems".

Monotonic basin hopping was developed by Leary [2] as a special case of the original basin hopping algorithm by Wales and Doye [1]. Algorithmically, it relies on the user specifying two probability distributions; one for the global search at Level-1 and another for local search at Level-2 as depicted in Fig. 1. Traditionally, the distributions have been chosen as uniform distributions with the global one having support over the full bounded parameter space and the local one having a user-tuned restricted width centered on the evolving sequence of algorithm iterates. Englander and Englander [17] provided a study of the efficiency and robustness of MBH with various local probability distributions when applied to a multiple gravity assist low-thrust (MGA-LT) Earth-Earth-Saturn-Uranus interplanetary transfer with variable specific impulse (VSI) propulsion and Sims-Flanagan transcription [18]. In particular, they looked at the performance when Pareto, Cauchy, Gaussian, and Uniform probability distributions are used and further examined the robustness to changes in each distribution's parameters. Numerical results showed that sampling from a Pareto distribution delivered the most robust performance, which the authors hypothesize is due to the heavy tailed characteristic of the distribution. Englander and Englander [19] have further explored the topology of the solution space for similar problems by analyzing the role of equality constraints on the sequence of iterates produced by a numerical solver and developing an adaptive hop probability distribution for MBH with application to an MGA-DSM mission [20].

One can view the work of Ampatzis and Izzo [21], and Cassioli et al. [22] as learning acceptance/rejection rules for samples from a globally supported uniform distribution at Level-1 in Fig. 1. Equivalently, one can view this as learning a uniform distribution with restricted support. Ampatzis and Izzo [21] achieve this by learning a two-layer artificial neural network (ANN) with 25 nodes for each layer to be a surrogate to the objective function of an MGA-DSM problem. Samples from a uniform distribution with poor evaluation from the surrogate are then rejected; improving numerical efficiency by avoiding costly evaluation for samples that are unlikely to yield high quality results from a gradient based numerical solver. Cassioli et al. [22] applies a nonlinear Support Vector Machine (SVM) classifier based on a Gaussian kernel to develop an acceptance/rejection criteria. The authors applied their approach to an MGA trajectory design problem and compared efficiency with that of MBH. In the work of Zhu and Luo [23], two multilayer perceptron ANNs with varying number of layers (2 to 5) and nodes (16 to 128) are trained to classify whether an unperturbed two-body low-thrust transfer is feasible in a fixed time and predict what the fuel consumption might be. The use of a classifier in their work is similar to the learning of the acceptance/rejection rules of Ampatzis and Izzo [21], and Cassioli et al. [22].

One of the key insights from the field of global optimization and the aforementioned studies in space flight is the fact that convergence of a global optimization algorithm requires global information (see for example the work of Stephens and Baritompa [24]). By extension, efficient global search for high quality solutions requires global information (see for example the chapter by Schoen [25]). In the case of the global search problem, the global information that is most pertinent is the topology of the basins and funnel structures, should they meaningfully exist. As the relative volume of the basins of attraction decrease due to increasing dimensionality of the searchable space, and as the geometry of the funnel structures become more sparse and less connected, the global search for high quality solutions becomes less likely to succeed with naive uniform sampling. Hence, global information about the geometry and topology should be used to further restrict sampling at Level-1 of Fig. 1. This paper demonstrates a framework that targets this concept by using generative machine learning to restrict sampling from a conditional probability distribution learned from prior solutions of similar problems.

## II. Mathematical Problem Formulation

When the global search process is over a parameterized set of problems, for example that arise in preliminary spacecraft trajectory mission design when system parameters, objectives, requirements, or constraints change, then a set of global search problems must be solved. Intuitively, if the topology of solutions of the global search problem change in a known manner (e.g., continuously) with respect to changes in the parameter variations, then a designer should make use of prior known solutions to accelerate the search for the newly parameterized problem. In this paper, we focus on exactly this problem. In particular, we propose a data-driven approach to learn a conditional probability distribution that can be used as the first-level algorithm in the global search process. To clearly articulate and justify the construction of our framework and learned distribution, as well as understand our final results, we must first describe how solutions of the lowest-level change as the parameterized global search problem and all choices at the lowest-level vary. This is the aim of the current section, which starts with the mathematical definition of a parameterized continuous-time optimal control problem (Sec. II.A) and eventually leads to the definition of local neighborhoods of high quality numerical solutions (Sec. II.D), upon which we can define properties for defining the conditional probability distribution that we seek to learn (Sec. II.E).

## A. Parameterized Optimal Control Problem

Our ultimate goal is to solve for a set of solutions to the collection of parameterized optimal control problems, which we denote by $\{O_\alpha\}$, and are given by Eqs. (1), (2), (3), and (4). In particular, the cost function $J(u; \alpha)$ to be minimized is parameterized by $\alpha$ and given by,

$$\min \left\{ J(u; \alpha) \equiv \phi(\xi(t_f), t_f; \alpha) + \int_{t_0}^{t_f} \mathcal{L}(\xi_s, u_s, s; \alpha) ds \,\middle|\, u \in \mathcal{U}, \text{ and Eqs. (2), (3), (4) are satisfied} \right\}, \quad (1)$$

where $\mathcal{U}$ is an admissible control set, $\phi$ is a terminal cost, and $\mathcal{L}$ is the running cost. The process $(\xi_s)_{s \in [t_0, t_f]}$ represents the state of the spacecraft and satisfies the ordinary differential equation (i.e., dynamical constraint),

$$\xi_t = \xi_0 + \int_{t_0}^{t} f(\xi_s, s; \alpha) ds + \int_{t_0}^{t} g(\xi_s, u_s, s; \alpha) ds, \quad \forall t \in [t_0, t_f], \quad (2)$$

with initial and terminal boundary conditions,

$$\xi_0 \equiv \xi(t_0) \in \Xi_0 \subseteq \mathbb{R}^m \quad \text{and} \quad \xi(t_f) \in \Xi_f \subseteq \mathbb{R}^m, \quad (3)$$

and a set of path constraints,

$$\psi_k(\xi_s, s; \alpha) \leq 0, \quad \forall s \in [t_0, t_f], \quad \forall k \in \mathcal{K}. \quad (4)$$

The vector field given by $g$ in Eq. (2) represents the perturbed forces due to control and $f$ is the vector field for any natural uncontrolled dynamics. The initial time $t_0$ and final time $t_f$ appear fixed in the equations, but could also be considered control parameters. The index set $\mathcal{K}$ is for the collection of path constraints. The boundary conditions $\Xi_0$ and $\Xi_f$ are assumed to be smooth submanifolds, which may not a have a global parametric description. For brevity we do not explicitly write the possible dependency of the parameter $\alpha$ on the admissible control set or boundary conditions, but it is reasonable and possible. We assume that all the coefficients in the aforementioned equations are nice (for example, that $\phi, \mathcal{L}, f, g$ and $\psi$ are smooth functions and that $\mathcal{U}$ is a collection of Lebesgue integrable functions), from which the standard theory for the existence of solutions may follow (see for example Bryson and Ho [26] or Liberzon [27]), and that solutions do indeed exist.

## B. Control Transcription

Rarely is it possible to solve a general optimal control problem analytically. It is therefore common to convert the continuous-time optimal control problem into a parameter optimization problem, typically a nonlinear program (NLP), by way of a control transcription. Doing so directly on Eqs. (1), (2), (3), and (4) results in the direct optimal control approach as opposed to the indirect approach that derives the first order optimality conditions and then transcribes the boundary value problem into an NLP. In either case, a choice of a control transcription is ultimately made, though less variation in the transcription is possibly for the indirect approach. We point the reader to some of the classical and common references on indirect [8, 26, 28], and direct methods [29–34], and the references therein for further details on various transcriptions. In this paper, we make use of a direct optimal control approach, using a forward-backward shooting control transcription to formulate an NLP (c.f., Sec. III.B). The NLP that is formed will be dependent on the original parameterized optimal control problem $O_\alpha$ as well as the control transcription. To make this clear, we let $h \in \mathbb{H}$ denote a control transcription from a collection of possible transcriptions given by the set $\mathbb{H}$. Our NLP then explicitly shows dependence on $\alpha$ and $h$ when written in the following form,

$$\min_{x \in \mathcal{X}^h} J(x; \alpha, h),$$
$$\text{subject to} \quad c_i(x; \alpha, h) = 0, \quad \forall i \in \mathcal{E}, \quad (5)$$
$$c_i(x; \alpha, h) \leq 0, \quad \forall i \in \mathcal{I}.$$

In Eq. (5) the objective function $J = J(\cdot; \alpha, h)$ and constraints $(c_i)_{i \in \mathcal{E} \cup \mathcal{I}}$ are assumed to be nonlinear functions and with sufficient regularity as stated in Sec. II.A. The objective function and constraints are functions of the state $x \in \mathcal{X}^h$, where $\mathcal{X}^h$ is a discretization of the admissible control set $\mathcal{U}$ given in Eq. (1) based on the transcription $h$. The sets $\mathcal{E}$ and $\mathcal{I}$ are index sets for equality and inequality constraints respectively. We denote the set of parameterized NLP problems as $\{\mathcal{P}_{\alpha,h}\}$.

## C. Solution via Numerical Solver

The solution of an NLP problem $\mathcal{P}_{\alpha,h}$ is by choice of a numerical solver (NS). We let the set of possible numerical solvers and their algorithm parameterizations by defined by the collection $\{\pi_\gamma\}$ indexed by $\gamma$. Assuming the solver to be gradient based, because the functions in Eq. (5) are differentiable, an initial (primal) state $x_0 \in \mathcal{X}^h$ is required for solution. NLP solvers often verify optimality of a proposed state $x$ by checking whether the Karush-Kuhn-Tucker (KKT) conditions [35, 36] hold. Therefore it is actually required that an initial state for the dual variables $\lambda_0 \in \Lambda$ is also provided*. Hence letting $(x_0, \lambda_0) \in \mathcal{X}^h \times \Lambda \equiv \mathcal{U}^h$ be an initial guess for an NLP $\mathcal{P}_{\alpha,h}$, an NS $\pi_\gamma$ will produce a sequence $(x_k, \lambda_k)$ of iterates that ideally converge to a local extremum $(x^*, \lambda^*)$ for a well-posed NLP. Therefore we think of an NS as a mapping from the space of possible NLPs and initial primal-dual states to an infinite sequence of primal-duals,

$$\pi_\gamma : (\mathcal{P}_{\alpha,h}) \times \mathcal{U}^h \longrightarrow \prod_{k=0}^{\infty} \mathcal{U}_k^h. \tag{6}$$

### D. High Quality and Diverse Extremum Solutions

*1. High Quality Solutions*

Our goal in pursuing global search during the preliminary mission design phase is to generate a large collection of high quality solutions that are also diverse. To make these ideas clear, we now fix the control transcription $h$ and the numerical solver $\pi_\gamma$, and define the subset $\mathcal{A}_\alpha \subset \mathcal{U}^h$ as the limiting points under the action of $\pi_\gamma$ that also satisfy the KKT conditions for $\mathcal{P}_{\alpha,h}$ within some tolerance, which may be included in $\gamma$. In particular, $\mathcal{A}_\alpha$ are the local extremum of $\mathcal{P}_{\alpha,h}$ under the action of $\pi_\gamma$. It is natural to define a solution as high quality if it belongs to $\mathcal{A}_\alpha$ and has an objective value better than some threshold. Therefore we define,

$$\mathcal{A}_{\alpha,\beta} \equiv \{(x, \lambda) \in \mathcal{A}_\alpha \mid J(x; \alpha, h) \leq \beta \in \mathbb{R}\}, \tag{7}$$

where $\beta$ serves as the threshold. As $\beta \to \infty$ we uncover the entire extremum set $\mathcal{A}_{\alpha,\beta} \to \mathcal{A}_\alpha$, and if $\beta < \min_{(x,\lambda) \in \mathcal{A}_\alpha} J(x; \alpha, h)$ then $\mathcal{A}_{\alpha,\beta} = \emptyset$.

When considering the design of intermediate-level algorithms for the global search problem, it is useful to consider multiple threshold levels $(\beta_j)_{j=0}^N$, where the collection $(\beta_j)$ has an increasing order $\beta_0 < \beta_1 < \cdots < \beta_N$, and $\mathcal{A}_{\alpha,\beta_j}$ is defined iteratively,

$$\mathcal{A}_{\alpha,\beta_0} \equiv \{(x, \lambda) \in \mathcal{A}_\alpha \mid J(x; \alpha, h) \leq \beta_0 \in \mathbb{R}\}, \tag{8}$$

$$\mathcal{A}_{\alpha,\beta_j} \equiv \{(x, \lambda) \in \mathcal{A}_\alpha \mid J(x; \alpha, h) \leq \beta_j \in \mathbb{R}\} \setminus \bigcup_{k=0}^{j-1} \mathcal{A}_{\alpha,\beta_k}, \quad j = 1, \ldots, N. \tag{9}$$

A directed graph may then be defined between the sets $(\mathcal{A}_{\alpha,\beta_j})$ to mathematically describe the likelihood of moving from some state $u \in \mathcal{A}_{\alpha,\beta_j}$ to a state $v \in \mathcal{A}_{\alpha,\beta_k}$ where $k < j$. Since it is beyond the scope of this paper to develop the intermediate-level algorithms, we say no more regarding this structure, and proceed with our analysis using a single threshold level $\beta$. The interested reader should consult the original clustering techniques and multilevel linkage ideas of

---

*Solvers typically have a default for the dual variables (e.g., a zero vector) when they are not provided

Becker and Lago [37], Törn [38], Rinnooy Kan and Timmer [39, 40] and the follow-on simple and random linkage methods of Locatelli and Schoen [41, 42] for further information on these perspectives.

*2. Diverse Solutions*

The meaning of diverse (or qualitatively different) solutions is a reflection of the fact that a mission designer is unable to completely model all desirable aspects of a mission solution into the original optimal control problem a priori. This may be due to a conflicting requirement of designing under high cadence, which promotes flexible lower fidelity models for quick solution and iteration, but also the reflection of reality that the design will change through iteration (i.e., the $\alpha$ parameters) and the designer simply does not know the design path a priori. Therefore one task of a mission designer is to review high quality solutions from $\mathcal{A}_{\alpha,\beta}$ and down select to a set with diverse solutions; that is, define a new subset $\overline{\mathcal{A}} \subseteq \mathcal{A}_{\alpha,\beta}$ such that each pair of solutions $u, v \in \overline{\mathcal{A}}$ are "qualitatively" different in some meaningful way (e.g., the solution corresponding to $u$ may traverse through a radiation belt, but $v$ avoids this region; or $u$ requires more revolutions than $v$ of a central body). The mission designer therefore has a "hidden" diversity "metric", which we denote as $\hat{d} : \mathcal{U}^h \times \mathcal{U}^h \to [0, \infty)$ that they use subjectively to construct $\overline{\mathcal{A}}$ from $\mathcal{A}_{\alpha,\beta}$. We assume that this diversity "metric" has the following properties,

$$\hat{d}(x, y) = \hat{d}(y, x) > 0 \quad \text{if} \quad x \neq y, \tag{10}$$

$$\hat{d}(x, x) = 0. \tag{11}$$

It is not a proper metric since we do not require the triangle inequality property. Assuming a diversity metric $\hat{d}$ exists, a mission designer can construct $\overline{\mathcal{A}}$ as the solution to the optimization problem,

$$\max_{\overline{\mathcal{A}} \subseteq \mathcal{A}_{\alpha,\beta}} |\overline{\mathcal{A}}|, \quad \text{subject to} \quad \min_{z \in \overline{\mathcal{A}}} \hat{d}\left(z, \overline{\mathcal{A}} \setminus \{z\}\right) > \eta > \eta^- > 0, \tag{12}$$

where $\eta > 0$ is a positive threshold set by the designer, and must be larger than,

$$\eta^- \equiv \min_{z \in \mathcal{A}_{\alpha,\beta}} \hat{d}\left(z, \mathcal{A}_{\alpha,\beta} \setminus \{z\}\right), \tag{13}$$

to be meaningful (we assume that $\mathcal{A}_{\alpha,\beta}$ is finite). The operator $|\cdot|$ in Eq. (12) gives the cardinality of the set. Therefore, the solution of Eq. (12) provides the largest collection of high quality solutions that are at least $\eta$-diverse.

Because we do not usually have a definition of $\hat{d}$ in practice, a surrogate can be convenient. A simple choice to make in this regard is to use any true metric (e.g., the standard Euclidean metric) defined on the space $\mathcal{U}^h \subset \mathbb{R}^n$ where $n = n(\alpha, h)$ is dependent on the parameterized optimal control problem $O_\alpha$ and the transcription $h$. This would appear to be a reasonable choice, since at least locally for any $u \in \mathcal{A}_{\alpha,\beta}$ one can expect, although not necessary, small variations in the solution state ($\xi_s$) for small changes in the control. In this paper, we take an applied approach of separating the priorities of finding a collection of high quality solutions and the further task of finding the maximal set of diverse solutions within this collection. This is accomplished by setting the $\beta$-threshold to be generous enough such that after finding $\mathcal{A}_{\alpha,\beta}$, we can manually inspect elements of the set and verify that diversity exists.

### E. A Conditional Probability Distribution for Efficient Global Search

With the preliminary definitions now complete, assume again that we have a fixed $h, \pi_\gamma$ and $\beta$ threshold, and we aim to learn a conditional probability distribution for the parameterized family $\{\mathcal{A}_{\alpha,\beta}\}_\alpha$. Given an $\alpha$ that we have not yet solved for, we would like to be able to quickly solve for the set $\mathcal{A}_{\alpha,\beta}$. There is a practical tradeoff that must be made in learning the conditional probability distribution $p = p(\cdot|\alpha)$. If the distribution is defined only with support on $\mathcal{A}_{\alpha,\beta}$, then initializing our NS $\pi_\gamma$ with $u \sim p(\cdot|\alpha)$ will result in immediate convergence. Hence a distribution with this support property is the ideal case for a quick warm-start. But if the set $\mathcal{A}_{\alpha,\beta}$ is finite, then it has Lebesgue measure zero, and the distribution $p$ will most likely vary in a complex manner as $\alpha$ varies, which implies requiring a very large amount of data to learn the conditional distribution; which may not be possible. At the other extreme, setting $p$ equal to the uniform distribution with support on all of $\mathcal{U}^h$ and for all $\alpha$, requires no data to learn, but provides no warm-start solution time improvements over the naive sampling from a uniform distribution (as we have trivially chosen it to be exactly the uniform distribution).

To define a conditional distribution that provides quick warm-start capability (i.e., few iterations of $\pi_\gamma$ for convergence to a local extremum) and hopefully requires a moderate amount of data, we introduce the idea of local neighborhoods. The $k$-local neighborhood of $\mathcal{A}_{\alpha,\beta}$ under $\pi_\gamma$ is defined as,

$$^k\mathcal{N}_{\alpha,\beta} \equiv \{z \in \mathcal{U}^h \mid \pi_\gamma^k(z) \in \mathcal{A}_{\alpha,\beta}\}, \tag{14}$$

where $\pi_\gamma^k(\cdot)$ indicates $k$-iterations of the solver (e.g., typically major iterations of an NLP solver). The zero neighborhood is such that $^0\mathcal{N}_{\alpha,\beta} = \mathcal{A}_{\alpha,\beta}$, and when $k \to \infty$ we have the local basin of attraction for $\mathcal{A}_{\alpha,\beta}$ under $\pi_\gamma$ (i.e., all points that eventually converge). The parameter $k$ in the above definition allows us to trade complexity in the definition of the conditional distribution with speed of convergence when samples from the conditional distribution are used in the warm-start of $\pi_\gamma$. Notice that the $k$-local neighborhood is equivalent to a union over the $k$-local neighborhoods of the elements of $\mathcal{A}_{\alpha,\beta}$:

$$^k\mathcal{N}_{\alpha,\beta} = \bigcup_{u \in \mathcal{A}_{\alpha,\beta}} {}^k\mathcal{N}_{\alpha,\beta}(u) \equiv \bigcup_{u \in \mathcal{A}_{\alpha,\beta}} \{z \in \mathcal{U}^h \mid \pi_\gamma^k(z) = u\}. \tag{15}$$

In particular, the sets $^k\mathcal{N}_{\alpha,\beta}(u)$ form a partition of $^k\mathcal{N}_{\alpha,\beta}$ for every $k$.

Let us denote $\mu_L$ as Lebesgue measure on $\mathcal{U}^h$ and assume that $\mu_L({}^k\mathcal{N}_{\alpha,\beta}(u)) > 0$ for all $k > 0$ and for all $u \in \mathcal{A}_{\alpha,\beta}$. The $k$-local neighborhoods generate an increasing nested set under the ordering given by $\mu_L$,

$$\mathcal{A}_{\alpha,\beta} \subseteq {}^1\mathcal{N}_{\alpha,\beta} \subseteq \cdots \subseteq {}^k\mathcal{N}_{\alpha,\beta} \subseteq \cdots, \tag{16}$$

and the same holds for each $k$-local neighborhood of $u \in \mathcal{A}_{\alpha,\beta}$.

One could similarly define $t$-local neighborhoods, where $t \geq 0$ is the time to converge. From an application perspective, this is more relevant, but requires that we account for the software implementation of $\pi_\gamma$ and the machine architecture that the implementation is running on. We use this notion of local neighborhoods along with the $k$-local neighborhood definition in Sec. V.B.1 when data is generated for learning the conditional distribution of our test problem.

Now that we have a definition of $k$-local (or $t$-local) neighborhoods to $\mathcal{A}_{\alpha,\beta}$, we can return to our goal of defining a conditional distribution $p$ with good warm-starting capability and requiring moderate data to learn $p$. In particular we would like to find a conditional distribution $p$ that smoothly approximates the following weighted Dirac distribution,

$$p(\cdot|\alpha) \approx \sum_{u \in \mathcal{A}_{\alpha,\beta}} \mu_L({}^k\mathcal{N}_{\alpha,\beta}(u))\delta_u(\cdot). \tag{17}$$

The weight of each Dirac distribution is the measure of its neighborhood, $\mu_L({}^k\mathcal{N}_{\alpha,\beta}(u))$, which can be approximated by a Monte Carlo sampling of $\mathcal{U}^h$:

$$\mu_L({}^k\mathcal{N}_{\alpha,\beta}(u)) \approx \frac{1}{N}\left|\{u_i \in {}^k\mathcal{N}_{\alpha,\beta}(u)\}\right|, \quad \left(u_i \sim U(\mathcal{U}^h)\right)_{i=1}^N. \tag{18}$$

In Eq. (18), $N$ is the number of Monte Carlo samples and $U(\mathcal{U}^h)$ is a uniform distribution with support on the entire control space $\mathcal{U}^h$.

The form of $p$ could be further improved by taking into account its definition on a collection of $k$-local (or $t$-local) neighborhoods and the diversity in solutions. By accounting for a collection of local neighborhoods, the regularity of $p$ can be further constrained. And by imposing diversity of solutions, support of $p$ on some local neighborhoods can be removed, which has the effect of more efficient online sampling when warm-starting (see for instance the literature in Sec. II.D.1). Attempts to further generalize the construction of $p$ in these directions are current research of the authors and beyond the scope of this paper. Our goal in this paper is simply to set an appropriate $\beta$, $k$ and $t$ that capture the funnels of high quality solutions with moderate convergence time under $\pi_\gamma$, upon which an appropriately tuned basin hopping algorithm could be applied to better explore the funnel and find the funnel global minima. We analyze the collection of discovered solutions after numerical solution to show that the ideas of Sec. II.D do approximately hold true;

that the standard Euclidean metric on $\mathcal{U}^h$ can be used as a proxy to a more refined diversity metric. Such an analysis is shown in Sec. III.B for the main test problem of this paper; solving of low-thrust circular restricted three-body transfers. In Sec. IV we detail how generative machine learning is then used to learn conditional distributions. Test results are shown and explained in Sec. V.

## III. Test Problems

Two test problems are described in this section. The first problem is a classical low dimensional problem from the global optimization literature, which we use to elucidate several core ideas of the paper; that a family of parameterized optimization problems with clustering type behavior in local optimal solutions can be learned using the generative machine learning framework to be described in Sec. IV. This problem does not require the full capability of the framework described in Sec. IV. The second problem is the LT circular restricted three-body problem of Sec. III.B and serves to illustrate how the topology of solutions to a continuous-time optimal control problem can possess clustering structure with many local basins of attraction; encapsulating the idea of multiple funnels. This problem also serves as the main benchmark to demonstrate the efficacy of the full framework.

### A. De Jong's Fifth Function

The De Jong's 5th function [43], also known as Shekel's foxholes function [44], is a classic parameterized test function used for benchmarking global optimization algorithms. The function can be parameterized such that it is characterized by a large number of local minima over a multi-modal landscape with sharp peaks and valleys, requiring careful algorithms to effectively navigate to the region of the global optimum. We use a version of De Jong's 5th function as a demonstrating example of our amortized global search (AmorGS) framework that learns the topological structure of local minima. The set of parameterized optimization problem $\mathcal{P}_{\alpha,h}$, corresponding to this objective function are defined as,

$$\min \left\{ J(\boldsymbol{x}; \alpha) = \left(0.002 + \sum_{i=1}^{n} \frac{1}{1 + (x_1 - \overline{A}(\alpha)_{1i})^6 + (x_2 - \overline{A}(\alpha)_{2i})^6}\right)^{-1} \Big| \boldsymbol{x} = (x_1, x_2) \in [-50, 50] \times [-50, 50] \right\}, \quad (19)$$

where the matrix $\overline{A}(\alpha)$, which is dependent on the parameter $\alpha$ and another matrix $A$, is defined as,

$$\overline{A} = \overline{A}(\alpha) \equiv \begin{bmatrix} \cos \alpha & -\sin \alpha \\ \sin \alpha & \cos \alpha \end{bmatrix} A, \quad A \equiv \begin{bmatrix} A_{11} & A_{12} & A_{13} & \ldots & A_{1n} \\ A_{21} & A_{22} & A_{23} & \ldots & A_{2n} \end{bmatrix}, \quad \text{and} \quad \alpha \in [0, \pi/2]. \quad (20)$$

Each column of $\overline{A}$ defines a local minima of the objective function $J$ with parameter $\alpha$. The conditional parameter $\alpha$ controls the rotation of the template solutions given by $A$. When $\alpha = 0$, the minima are simply the columns of the matrix $A$. An example of this function with 25 local minima is depicted in Fig. 2.

## B. Low Thrust Circular Restricted Three-Body Problem

The second test problem of this paper considers a minimum-fuel transfer of a LT CSI spacecraft in the Earth-Moon system. The spacecraft's motion is modeled using the Circular Restricted Three-Body Problem (CR3BP) dynamics, which describes the motion of a point (negligible) mass under the gravitational influence of two celestial bodies. The CR3BP serves as a valuable low-fidelity model for preliminary mission design, providing a good first-order approximation of the complex dynamics in multibody dynamical environments (e.g., the Earth-Moon cislunar space). Following the standard construction, we assume that the Earth and Moon follow circular orbits relative to their barycenter, and we write the equations of motion of the spacecraft in a synodic reference frame rotating at the same rotation rate as the Earth and Moon. Non-dimensionalization is carried out to obtain a suitable choice of units, which reduces the number of parameters in the problem to one, namely, the mass parameter $\mu = m_2/(m_1 + m_2)$, where $m_1$ is the mass of the primary (e.g., Earth) and $m_2 \leq m_1$ is the mass of the secondary (e.g., Moon). With this choice of units, the gravitational constant and the mean motion both become unity, leading to the following equations of motion,

$$\ddot{\boldsymbol{q}} = -2\hat{\boldsymbol{q}}_3 \times \dot{\boldsymbol{q}} + \nabla_{\boldsymbol{q}} \overline{U}(\boldsymbol{q}), \tag{21}$$

where $\boldsymbol{q} = (q_1, q_2, q_3) \in \mathbb{R}^3$ describes the position in the synodic reference frame and $\hat{\boldsymbol{q}}_3$ is the unit vector normal to the orbital plane of the primary and secondary bodies. The effective gravitational potential $\overline{U}$ is defined as,

$$\overline{U}(\boldsymbol{q}) \equiv \frac{1}{2}\left(q_1^2 + q_2^2\right) + \frac{1-\mu}{r_1(\boldsymbol{q})} + \frac{\mu}{r_2(\boldsymbol{q})}, \tag{22}$$

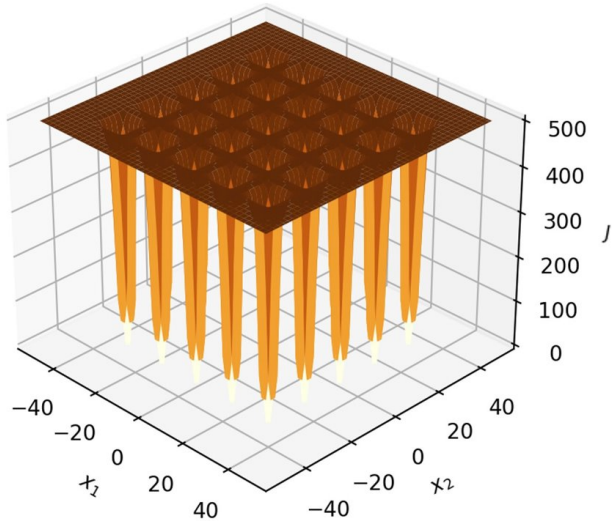

**Fig. 2  An example of De Jong's 5th function in $\mathbb{R}^2$ with 25 local minima.**

with $r_1$ and $r_2$ the distances between the spacecraft to the primary and the secondary,

$$r_1(\boldsymbol{q}) = \sqrt{(q_1 + \mu)^2 + q_2^2 + q_3^2}, \tag{23}$$

$$r_2(\boldsymbol{q}) = \sqrt{(q_1 - (1 - \mu))^2 + q_2^2 + q_3^2}. \tag{24}$$

For the Earth-Moon system, $\mu \approx 0.01215$, and the CR3BP model possess five relative equilibrium. The relative equilibrium that is collinear and between the Earth and Moon in the synodic frame is of interest in this paper, and is known as the $L_1$ Lagrange point, which possesses families of periodic spatial Halo orbits that have stable manifolds.

*1. Optimal Control Problem Definition*

As stated, we consider a minimum-fuel LT CSI transfer problem and therefore model the spacecraft mass $m$ as,

$$\dot{m} = -\alpha \frac{|u|}{I_{sp}\bar{g}}, \qquad |u| \leq T_{max}, \tag{25}$$

where $\bar{g} \approx 9.80665$ m/s$^2$ is the standard acceleration on Earth, $T_{max}$ is the maximum available thrust, and we recognize what will be the conditional parameter for families of solutions to this problem with the parameter $\alpha \in [0.1, 1]$. In Eq. (25), $u$ is the LT control, which then modifies the natural dynamics of Eq. (21) as follows,

$$\ddot{\boldsymbol{q}} = -2\hat{\boldsymbol{q}}_3 \times \dot{\boldsymbol{q}} + \nabla_{\boldsymbol{q}} \bar{U}(\boldsymbol{q}) + \alpha \frac{u}{m}. \tag{26}$$

The cost function to be minimized is,

$$J(u; \alpha) = -m(t_f) = -\int_{t_0}^{t_f} \dot{m}_s ds. \tag{27}$$

The constant specific impulse is $I_{sp} = 1000$ seconds, initial fuel mass is set equal to 700 kg, and dry mass equal to 300 kg. A maximum thrust of $T_{max} = 1$ N is used, which then gives a range of problems to be consider with $\alpha T_{max} \in [0.1, 1]$ N. Thus the initial maximum acceleration of the spacecraft ranges from $10^{-4}$ m/s$^2$ to $10^{-3}$ m/s$^2$. The trajectories of low-thrust spacecraft exhibit a high sensitivity to the maximum allowable thrust, with varying maximum thrust levels leading to qualitatively distinct trajectory patterns. Lastly, the boundary conditions for the spacecraft trajectory are chosen such that the spacecraft begins after a low-thrust spiral originating from a geostationary transfer orbit and terminates at a stable manifold arc of a Halo orbit to the $L_1$ Lagrange point.

*2. Control Transcription*

A forward-backward finite-burn low-thrust shooting transcription is used in this work to formulate a nonlinear program $\mathcal{P}_{\alpha,h}$ as in Eq. (5) with the decision variable $x \in \mathbb{R}^{3N+4}$, where $N$ will represent the number of control segments. With this forward-backward shooting transcription, the decision variable $x$ for the optimization problem is defined as follows,

$$x = (\tau_s, \tau_i, \tau_f, m_f, u_1, \ldots, u_N). \tag{28}$$

As a finite-burn low-thrust transcription, the spacecraft thrust vector remains constant across each segment characterized by the variable $u_i \in [0, 2\pi] \times [-\pi/2, \pi/2] \times [0, T_{\max}] \subseteq \mathbb{R}^3, i \in \{1, 2, \ldots, N\}$, representing the direction and magnitude of the thrust vector. In this paper the number of control segments is fixed at $N = 20$. The remaining four variables in Eq. (28) are the initial coast time $\tau_i$, the final coast time $\tau_f$, the shooting time $\tau_s$, and the final mass $m_f$ at the terminal boundary. The length of time of each control segment is $\Delta t = \tau_s/N$. Equality constraints are introduced for this problem to enforce continuity of the position, velocity, and mass of the spacecraft. As a forward-backward transcription, this occurs at the mid-point in the control sequence. In particular, the left boundary condition from the forward arc is,

$$\xi^-_{\tau_i+\tau_s/2} = \xi_0 + \int_{t_0}^{\tau_i} f(\xi_s, s)ds + \sum_{j=1}^{N/2} \int_{\tau_i+(j-1)\Delta t}^{\tau_i+j\Delta t} f(\xi_s, s) + g(\xi_s, u_j, s; \alpha)ds, \tag{29}$$

with $t_f = \tau_i + \tau_s + \tau_f$, the right boundary condition from the backward arc is,

$$\xi^+_{\tau_i+\tau_s/2} = \xi_{t_f} - \int_{t_f-\tau_f}^{t_f} f(\xi_s, s)ds - \sum_{j=1}^{N/2} \int_{t_f-\tau_f-j\Delta t}^{t_f-\tau_f-(j-1)\Delta t} f(\xi_s, s) + g(\xi_s, u_{N-j+1}, s; \alpha)ds, \tag{30}$$

yielding the equality constraints,

$$c(x) \equiv \xi^+_{\tau_i+\tau_s/2} - \xi^-_{\tau_i+\tau_s/2} \in \mathbb{R}^7, \tag{31}$$

with feasible solutions $x$ yielding $c(x) = 0$. The problem is solved using **pydylan**, which is the Python interface to the Dynamically Leveraged Automated (N) Multibody Trajectory Optimization solver (DyLAN) [45]. The solution of the NLP is done with the Sparse Nonlinear Optimization (SNOPT) software [46].

*3. Structure of Solutions*

To understand the topological structure of minima, we analyze a subset of the 300,000 solutions solved for in Sec. V.B.1. In particular, we consider just three cases for $\alpha$ in this section, $\alpha \in \{0.1, 0.3, 1.0\}$. We select 1,000 local optimal solutions via random uniform sampling for each of these different $\alpha$ scenarios and project their solutions into the space

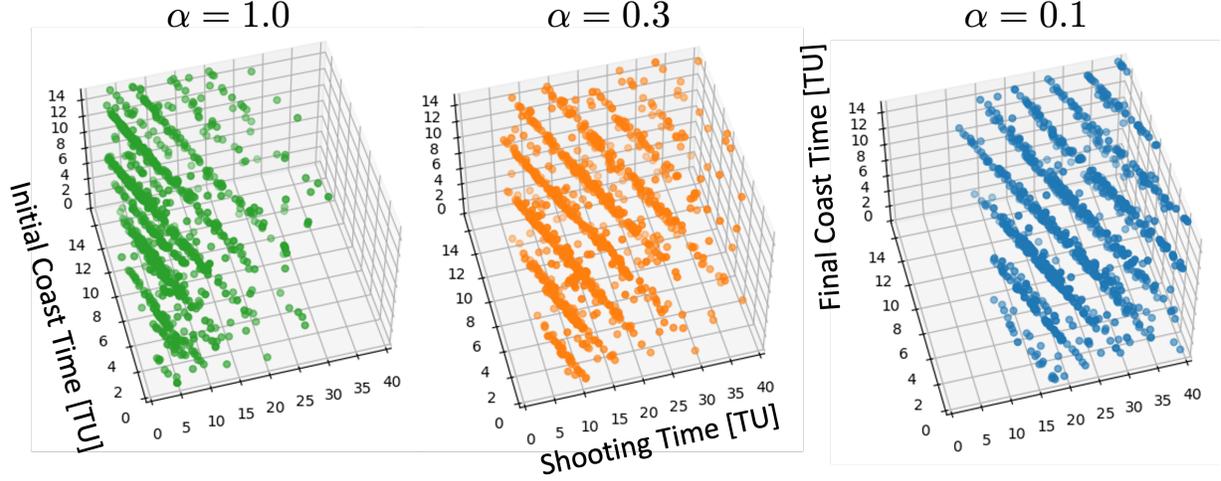

**Fig. 3** Projecting the time coordinates of solutions in the sets $\mathcal{A}_\alpha$ for $\alpha \in \{0.1, 0.3, 1.0\}$ reveals hyperplane structures. Decreasing $\alpha$ results in longer times-of-flight and translating hyperplanes.

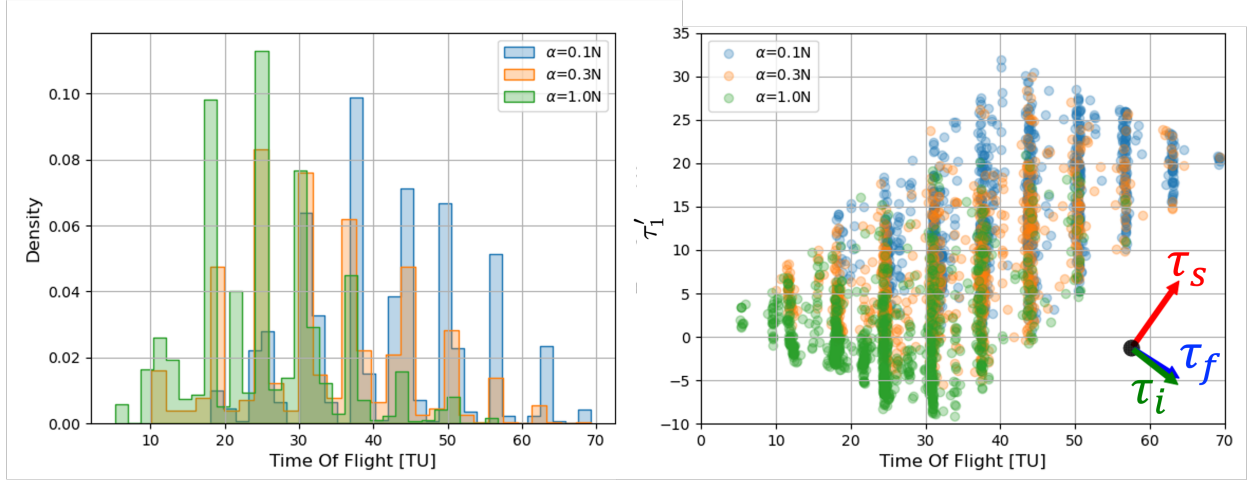

**Fig. 4** (Left) The histograms of the time-of-flight also shows a shift in the distribution as $\alpha$ changes; (Right) An appropriate transformation of hyperplanes helps to better reveal their structure.

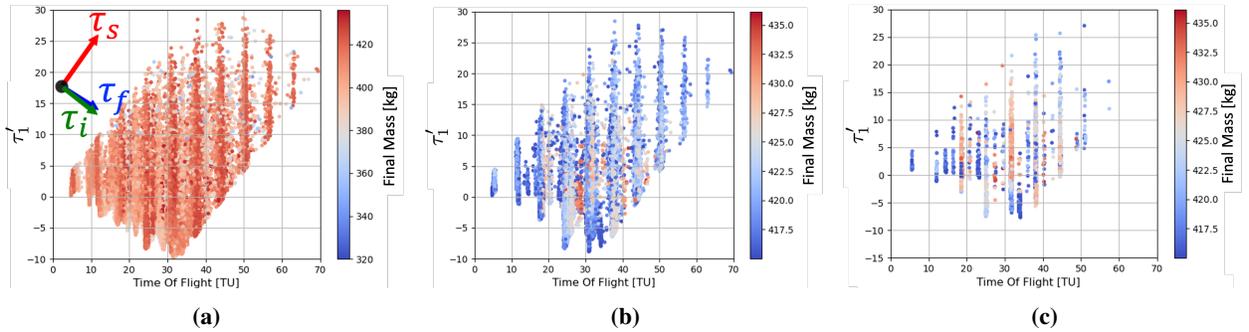

**Fig. 5** (a) Unfiltered data of $\mathcal{A}_{1.0}$, which contains **88,000 locally optimal solutions**, (b) Data filtered by a mass threshold of $\beta = 415$ kg (i.e., $\mathcal{A}_{1.0,415}$), (c) The $\mathcal{A}_{1.0,415}$ set further filtered by closeness to 20 hyperplanes (modes) shown in Fig. (4). All subfigures show data post-coordinate transformation to demonstrate the distinct hyperplane structure.
18

of time coordinates. Specifically, for each minima the projected coordinates are $\boldsymbol{\tau} = (\tau_s, \tau_i, \tau_f)$ and we let $\mathcal{F}_1$ denote this canonical time coordinate frame. Figure 3 shows the results of this reduced data set and in particular, the existence of hyperplane structures. The sum of the three time variables $\tau_s, \tau_i, \tau_f$, which corresponds to the total time-of-flight, remains approximately constant within each hyperplane. Figure 4 presents a histogram of this time-of-flight data. From this empirical distribution for each $\alpha$ we can identify hyperplanes as modes of the distribution. These figures begin to corroborate the hypothesis, as discussed in Sec. I, that space flight optimal control problems may generically contain clustering structure.

To better perform the in-depth analysis within and of adjacent hyperplanes, we perform coordinate transformations of the time data. The application of the coordinate transformations are shown in Figs. 4 and 5. Fixing an $\alpha$, let $T$ denote a mode of the $\alpha$-density as shown in Fig. 4. We define the hyperplane $H_T$ to this mode as,

$$H_T \equiv \{\boldsymbol{\tau} \mid |\|\boldsymbol{\tau}\|_1 - T| \leq \delta\}, \tag{32}$$

where $\delta$ is chosen to be 0.25 TU and $\|\cdot\|_1$ is the $l_1$-norm. To generate a coordinate transformation for the hyperplane $H_T$, three data points of $H_T$ are selected via uniform random sampling. Generically, the three data points span the hyperplane and a normal direction of the hyperplane is generated via the cross product of two in-plane orthonormal vectors. Denoting these orthonormal vectors as $(\hat{\tau}_1, \hat{\tau}_2, \hat{\tau}_3)$ and the reference frame they create as $\mathcal{F}_2$, the coordinate transformation $\mathcal{R}_{\mathcal{F}_1}^{\mathcal{F}_2} : \mathcal{F}_1 \to \mathcal{F}_2$ for $H_T$ is then,

$$\boldsymbol{\tau}' = \mathcal{R}_{\mathcal{F}_1}^{\mathcal{F}_2} \boldsymbol{\tau}. \tag{33}$$

The right subfigure of Fig. 4 shows the data of Fig. 3 after these transformations. For data not within the tolerance $\delta = 0.25$ TU of a given hyperplane, the coordinate transformation for the nearest hyperplane is used to transform.

We now restrict to the case of $\alpha = 1.0$, where we identify 20 modes and perform further analysis to understand the variation in the solution structure within each hyperplane and of adjacent ones. The features that we now discuss our qualitatively similar for $\alpha \neq 1.0$ and therefore we restrict our more in-depth analysis to this case. Fig. 5(a) shows the entire $\mathcal{A}_{1.0}$ data set of 88,000 minima in the transformed coordinates. Filtering out non-high quality solutions by setting our threshold to $\beta = 415$ kg yields Fig. 5(b)[†]. Figure 5(c) is the result of removing any data from Fig. 5(b) that is not within a $\delta = 0.25$ TU normal distance to a hyperplane and projecting the remaining data orthogonally onto the hyperplane.

Figure 6 shows the in-plane objective function structure of eight hyperplanes from Fig. 5(c). A contour plot is used on each hyperplane with a linear interpolation scheme using the `contourf` function in the `matplotlib` library [47].

---

[†]note that due to the cost function of Eq. (27), higher final mass solutions are better than lower ones

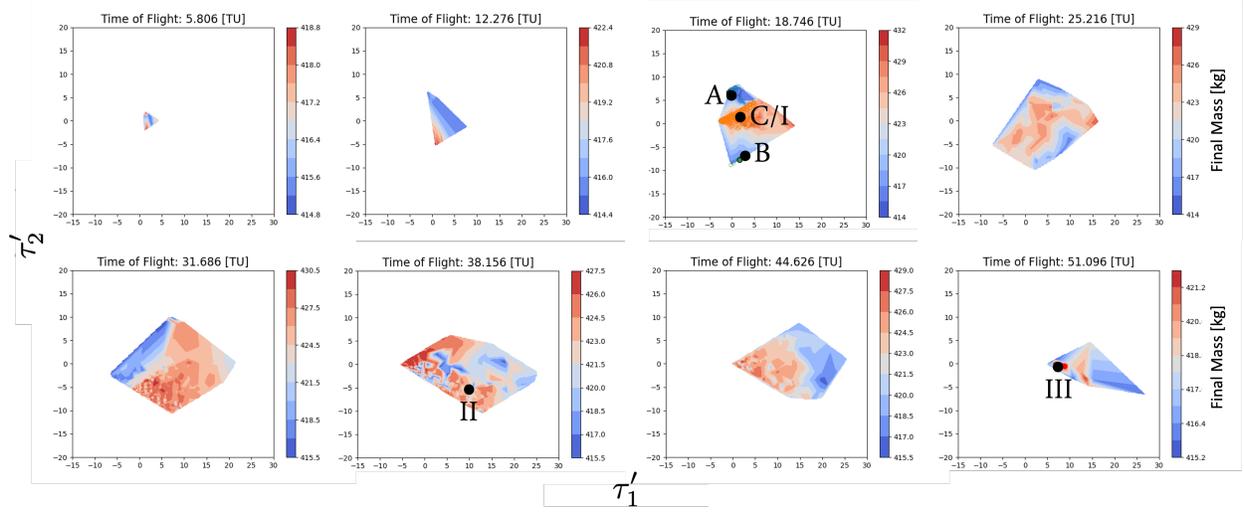

**Fig. 6    The distribution for the final mass of minima are shown across varying hyperplanes corresponding to different times-of-flight.**

Ten contour levels are shown in the figures. Fixing a given hyperplane, we see different clustering of basins of attraction resulting in possible funnel structure. The 18.746 TU slice shows what appears to be a single funnel to the center of the hyperplane in the $\tau_1'$ and $\tau_2'$ coordinates. The hyperplane for 25.216 TU depicts a possible funnel structure that spreads in several directions. The hyperplanes of 31.686, 38.156, and 44.626 TU all have more localized effects with what might be many local minima in single or double funnel structures. This analysis is of course limited, only showing two degrees of freedom, and should at least be compared with Fig. 5(c) that shows certain hyperplanes having minima with superior objective values over others isolated from each; that is, there existence some hyperplanes between those with the best minima, which may imply a multiple funnel structure.

*4. An Averaging Technique to Reveal Funnel and Multi-Modal Structure*

Figure 1 provides a nice visual of hypothesized structure for the objective function value of a multiple funnel multi-modal problem over a one-dimensional control problem. To realize a similar technique for the problem at hand, we apply moving averages (MA) over the time-of-flight for the full data set $\mathcal{A}_{1.0}$ with various sliding window sizes. This technique allows us to attempt to capture the topological properties of the solution landscape and reveal its multi-scale structure in an averaged sense.

The largest MA window sizes of 0.25, 0.5, and 1.0 TU are shown in Fig. 7. All of the data for $\mathcal{A}_{1.0}$ is shown as black markers. The MA trends show increasing frequency variation as the MA window sizes decrease. At the macro scale, they imply multiple funnels with increasing number of local minima traversing to the top[‡] of the funnels. In Fig. 8 we provide a closer look at the fine and intricate averaged structure centered on 24.5 TU. Assuming that the structure shown here is independent of the non-time control parameters, one can devise efficient global search algorithms that search with step sizes large enough to overcome the finest variability shown in the MA windows of 0.01 TU, yet doesn't overstep the important basins appearing at the 0.05 TU MA. Similarly, for the construction of a conditional probability distribution to be sampled at Level-1 in Fig. 1, an efficient distribution would have support on basins of approximately 0.05 TU in width.

---
[‡]note that here the objective function is the negative of that solved, hence as shown, the objective function would have been a maximum

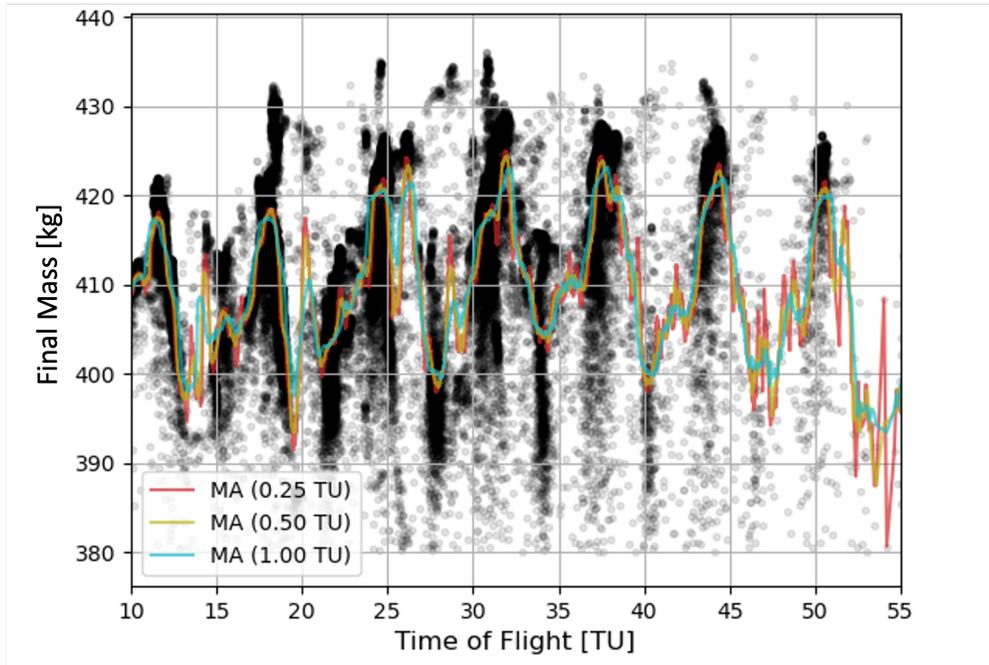

**Fig. 7** Moving average of the objective function over time-of-flight with different moving average windows, revealing macro-level funnel structure and multi-modality.

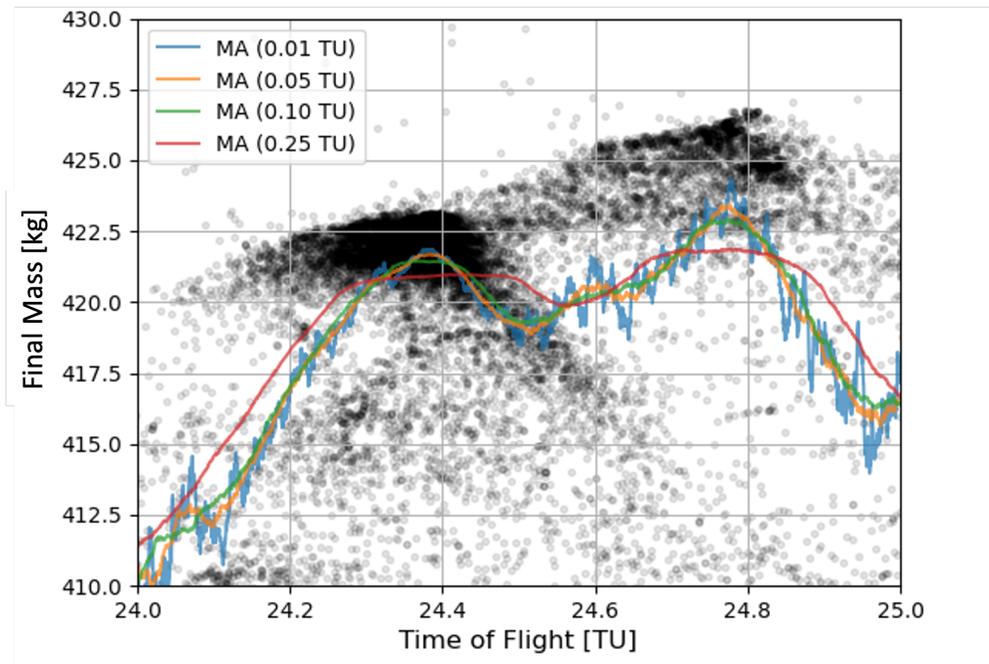

**Fig. 8** Moving average of the objective function over time-of-flight with narrow moving average windows, revealing micro-level funnel structure and multi-modality.

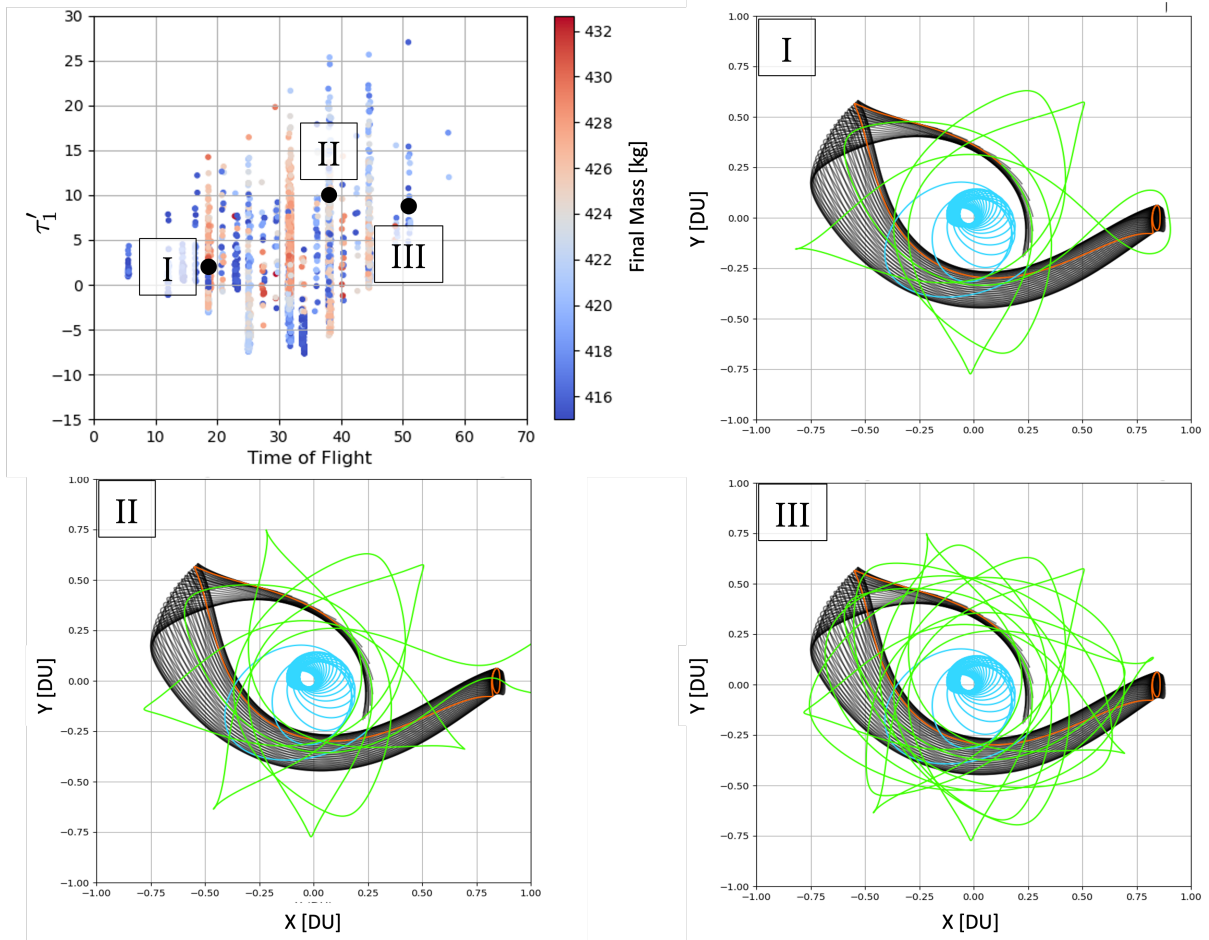

**Fig. 9  Three qualitatively different local minima for the minimum-fuel LT CR3BP transfer problem $\mathcal{P}_{1.0,h}$, sampled from hyperplanes with different times-of-flight. The Earth-Moon CR3BP distance unit is 1 DU = 384,400 km. Trajectories shown include: low-thrust spiral from geostationary orbit (in cyan), optimal control solution (in green), target manifold arc (in orange), and the stable invariant manifold (in black) of a Halo orbit around the $L_1$ Lagrange point.**

*5. Diversity across Hyperplanes (Macro Structure)*

In Section II.D.2, we emphasize the desire for qualitatively diverse trajectory solutions, which we assume here to loosely imply differences in the state solution $\xi_t(u_t)$ given a control $u_t$. In Fig. 9, we illustrate three solutions from $\mathcal{A}_{1.0}$ selected from three different hyperplanes, with times-of-flight equal to 18.746 TU (Solution I), 38.156 TU (Solution II), and 51.096 TU (Solution III) that display diverse qualitative solution behavior. The solutions in each hyperplane are labeled in Fig. 6 and are all high quality solutions with approximately equivalent delivered final mass. With regard to diversity, Solution I enters the Moon realm but remains on the near side, Solution II completes a lunar flyby and extends to the far side, and Solution III remains confined within the interior realm. Moreover, these solutions demonstrate varying numbers of revolutions, which correlate with the increasing times-of-flight.

## 6. Diversity within Hyperplanes (Micro Structure)

On the other hand, solutions belonging to the same hyperplane tend to be more similar to each other. In Fig. 10, we present three solutions from $\mathcal{A}_{1.0}$, all selected from the same hyperplane and thus having approximately identical times-of-flight equal to 18.746 TU. The location of these solutions, labeled as A, B, and C, are shown in Fig. 6. These trajectories are qualitatively similar, but exhibit slight variations. For instance, Solution I does not include any retrograde motion (i.e., 'loops'), whereas both Solution II and Solution III display this behavior. Since all solutions correspond to approximately the same time-of-flight, the number of revolutions remains comparable across them.

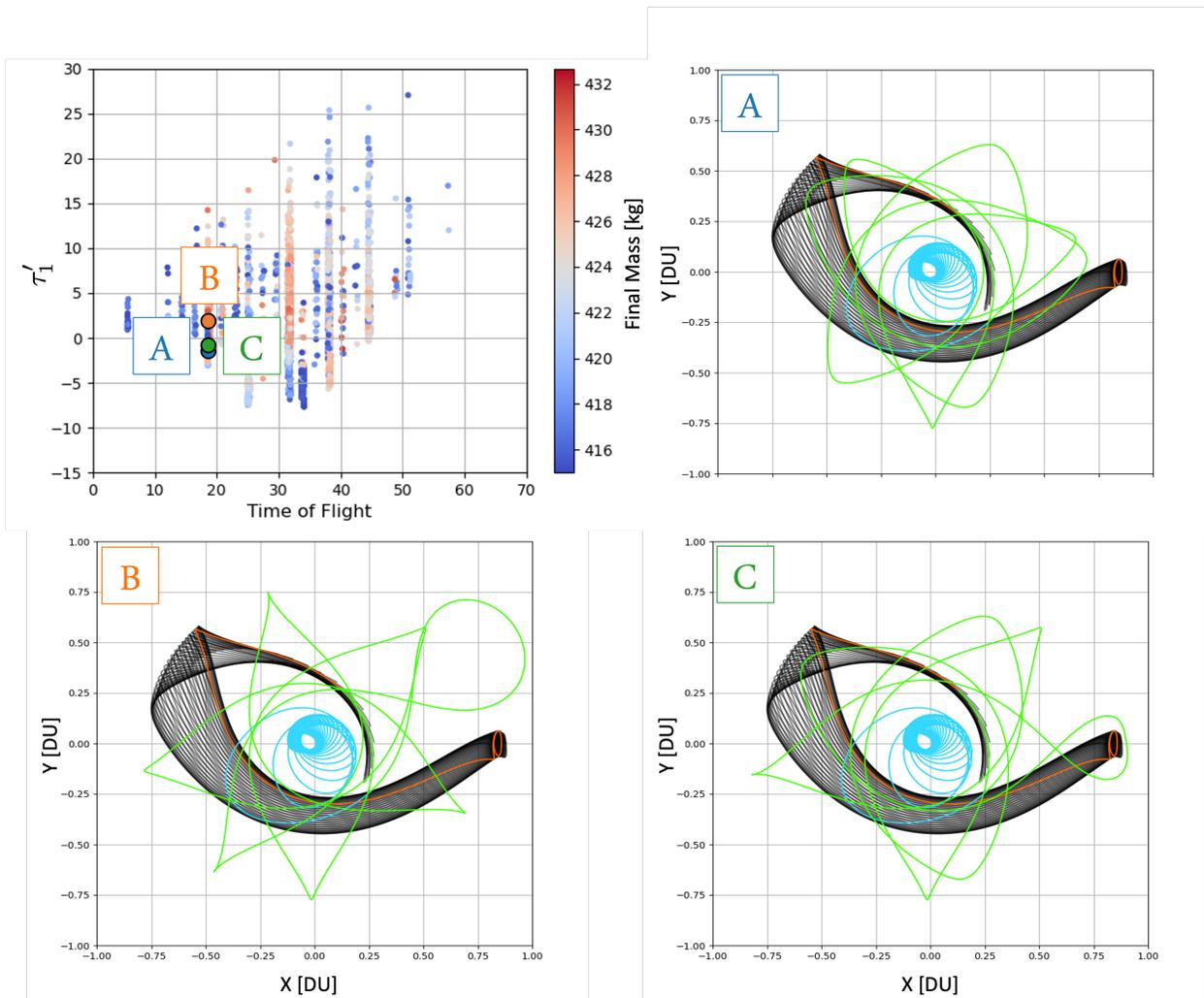

**Fig. 10** Three qualitatively similar local optimal trajectory solutions to the minimum-fuel LT CR3BP transfer problem $\mathcal{P}_{1.0,h}$, sampled from the same 18.746 TU hyperplane. See the caption of Fig. 9 for more details.

Prying more into the data of Fig. 6, we can visualize the raw solution data side-side with the contour plots to better understand if a given hyperplane has a single funnel or potentially multiple funnels. Figures 11 and 12 provide examples of this. Figure 11 shows up to three possibly distinct funnels for the hyperplane corresponding to an 18.746 TU time-of-flight, whereas Fig. 12 has a more complex structure for the 38.156 TU case.

In Fig. 11, the highest quality solutions exist in Funnel B. Funnel C is separated from Funnel B by a substantial distance in the $\tau_1', \tau_2'$ coordinates. Yet Funnel A may actually continue into what is denoted as Funnel B, implying that initial guesses for Funnel A and an appropriately tuned Level-2 algorithm can quickly traverse into the Funnel B location. If solutions in distinct funnels exhibit minor differences in their trajectories, as seen in Fig. 10 for the potential funnel areas of Fig. 11, then it could be advantageous for any sampler of a Level-1 algorithm to avoid sampling in what

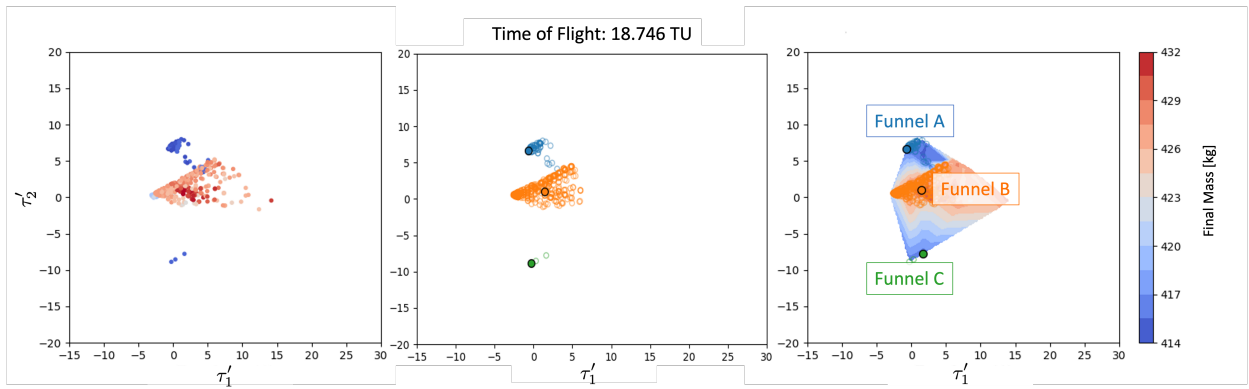

**Fig. 11  Potential funnel structure for the 18.746 TU hyperplane with raw data colorized by final mass (left), data categorized and colored by funnel (middle), and contour overlay of final mass (right).**

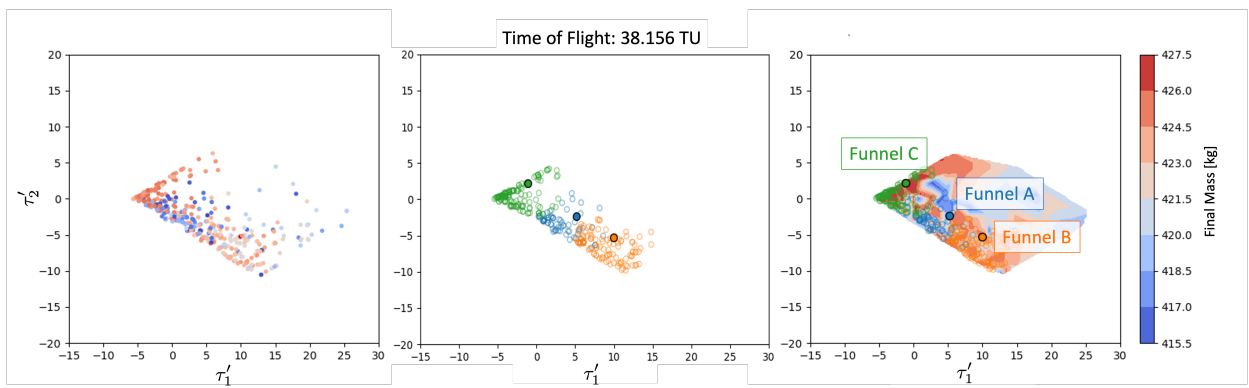

**Fig. 12  Potential funnel structure for the 38.156 TU hyperplane with raw data colorized by final mass (left), data categorized and colored by funnel (middle), and contour overlay of final mass (right).**

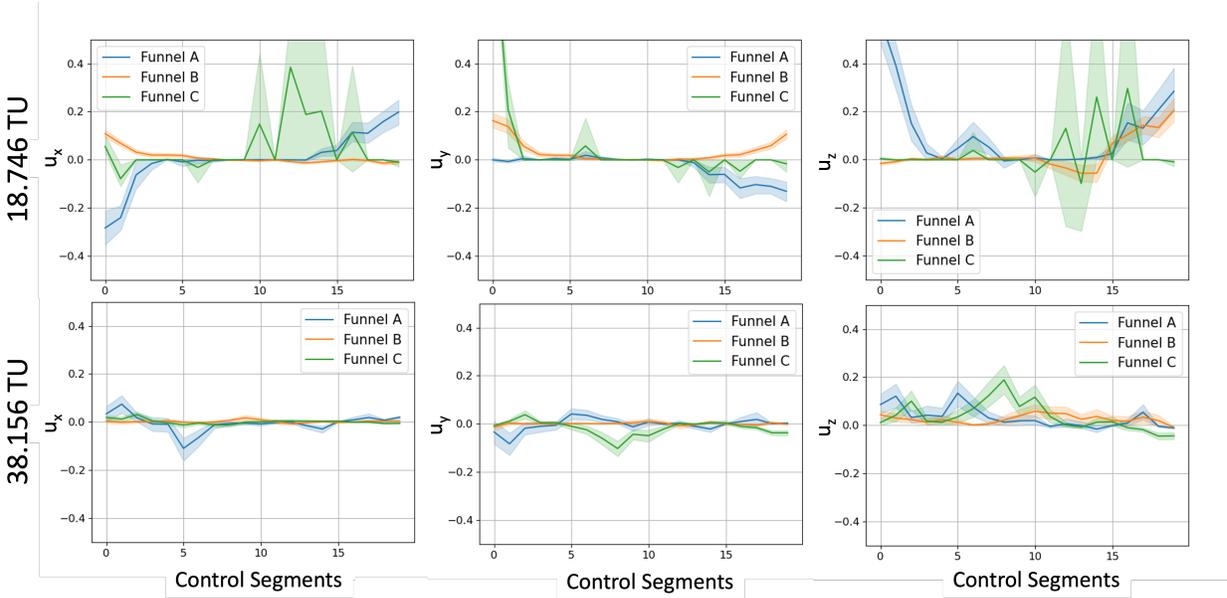

**Fig. 13** The mean and 95% confidence interval is shown for the thrust control components, in synodic coordinates, for the 18.746 TU and 38.146 TU hyperplanes across three funnels shown in Figs. 11 and 12.

is denoted as Funnel C of Fig. 11.

The same analysis for Fig. 12 is more complex and it must be remembered that these figures only show a subset of the full control parameter space. It is possible that Funnel B and Funnel C, as designated in Fig. 12 are distinct. The area denoted as Funnel A may in fact be appropriate continuations of either B or C, or possibly its own distinct funnel possessing solutions with objectives vales being slightly higher than the best seen in B or C. Still, the solutions in Funnel A, may be desirable for a mission designer if they provide qualitatively diverse solutions from those found in B or C.

## 7. Temporal Correlation

Thus far we have only investigated the solution space related to four of the control parameters in the transcription given in Eq. (29). The remaining parameters are associated with the thruster pointing direction and throttle ($u_1, \ldots, u_{20}$). As will be detailed in Sec. IV.D and the comparative results of Sec. V.B, it is important in the modeling architecture to account for the temporal correlation in these control segments. Doing so can drastically improve the efficiency of warm-starting the Level-3 numerical solver via a Level-1 sampler accounting for this structure. In Fig. 13 the mean and 95% confidence interval of the solutions in the basins identified in Figs. 11 and 12 are shown for each control direction over the number of control segments.

It is the solutions in Funnel B of Fig. 11 that have the least fuel consumption for the 18.746 TU case, even though these solutions have a larger average shooting time in comparison to solutions in Funnel A or C as shown in Table 1a. The fact that Funnel B solutions use less fuel is also seen in Fig. 13, where the control magnitude of Funnel B is most often close to zero magnitude in all components. The large variance in the Funnel C solutions is due to the few samples

**Table 1   Average time parameters across funnels presented in Figs. 11 and 12. The total coast time corresponds to the sum of initial and final coast times.**

| | Time of Flight = 18.746 TU | | | Time of Flight = 38.156 TU | |
| --- | --- | --- | --- | --- | --- |
| | Shooting Time [TU] | Total Coast Time [TU] | | Shooting Time [TU] | Total Coast Time [TU] |
| Funnel A | 6.479 | 12.066 | Funnel A | 15.252 | 22.879 |
| Funnel B | 6.791 | 11.777 | Funnel B | 19.960 | 18.124 |
| Funnel C | 6.648 | 11.890 | Funnel C | 10.671 | 27.554 |
| (a) | | | (b) | | |

**Table 2   Number of solutions in each funnel (see Figs. 11 and 12)**

| | Time of Flight [TU] | |
| --- | --- | --- |
| | 18.746 | 38.156 |
| Funnel A | 63 | 71 |
| Funnel B | 381 | 112 |
| Funnel C | 3 | 143 |

in that subdomain, as documented in Table 2.

The same analysis of the 38.146 TU case is richer. Here we see that the best solutions lie in Funnels B or C. Funnel B shows a nearly zero magnitude control authority throughout all control segments in Fig. 13, whereas the Funnel C shows significant deviation from zero magnitude, especially in the out-of-plane component that is necessary to transfer from the end of the low-thrust spiral from the GTO orbit to the Halo manifold arc. But we see from Table 1b that the Funnel C shooting time is also substantially less than solutions in B or A, and hence the integration of the larger control magnitudes is not as impactful on fuel usage. With substantial samples in each of the three subdomains, as shown in Table 2, the variances of the 38.156 TU case are much more consistent.

## IV. Amortized Global Search (AmorGS)

### A. AmorGS Framework

As stated in Sec. II.E, the AmorGS framework aims to learn a conditional probability distribution $p(\cdot|\alpha)$ for the collection of local optimal solution sets $\{\mathcal{A}_{\alpha,\beta}\}_\alpha$ in Eq. (7) to problems $\{\mathcal{P}_{\alpha,h}\}_\alpha$ of Eq. (5). Such a distribution is learned based on data generated from solved instances of $\{\mathcal{P}_{\alpha,h}\}_\alpha$, which may not include all $\alpha$ of interest. In the AmorGS framework, conditional generative machine learning (ML) is used to represent the condition distribution. In this work, we make use of a conditional variational autoencoder (CVAE) for this purpose. This was first proposed in Li et al. [48]. Other popular generative ML models include generative adversarial networks [49] and diffusion models [50–52], the latter of which has also been applied in the AmorGS framework recently [53–55]. There may be additional neural networks (NN) that are used to complete the generation of samples for the conditional distribution. In the main problem of this paper, given in Sec. V.B, a Long Short-Term Memory (LSTM) model is used as an additional NN.

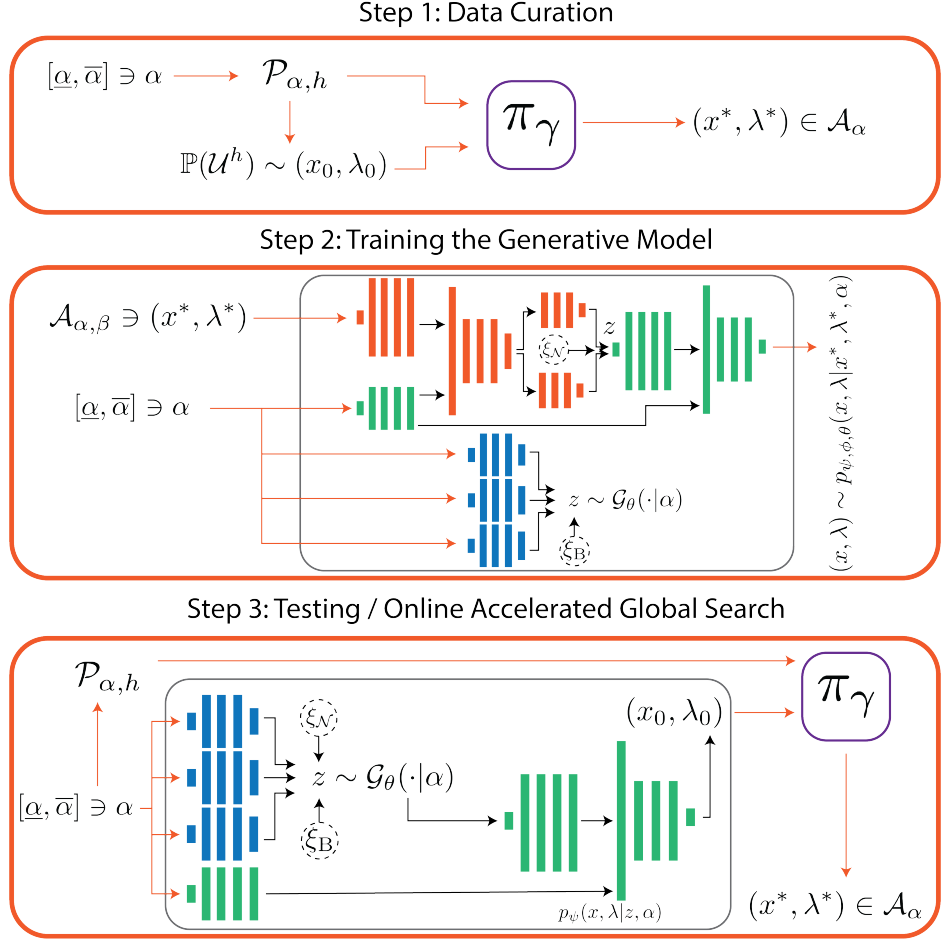

**Fig. 14** The workflow of the AmorGS framework to accelerate a global search over $\alpha$-parameterized optimization problems $\mathcal{P}_{\alpha,h}$ with numerical solver $\pi_\gamma$ by learning a parameterized conditional distribution $p_{\phi,\theta}(x,\lambda|\alpha)$ **given by a generative ML model.**

The main workflow of the AmorGS framework is shown in Fig. 14 and consists of three steps. In the first step, data is curated for instances of problems $\mathcal{P}_{\alpha,h}$ parameterized by $\alpha \in [\underline{\alpha}, \overline{\alpha}]$, where $\underline{\alpha}$ and $\overline{\alpha}$ are lower and upper bounds of a set interval of interest. The parameter $\alpha$ could be an element of a more general index set, but here we consider just an interval for simplicity. The data that is curated is due to the NS $\pi_\gamma$ acting on the problem $\mathcal{P}_{\alpha,h}$ with initial guess $(x_0, \lambda_0)$ and terminating with a solution $(x^*, \lambda^*) \in \mathcal{A}_\alpha$ that is a local extrema. In the second step, the curated data is filtered to retain high quality solutions specified by the threshold parameter $\beta$. This data is then used to train a generative model that is able to produce new initial guess samples $(x, \lambda)$, which should ideally be representative of the data set $\mathcal{A}_{\alpha,\beta}$. In the third and final step, the generative model takes as input $\alpha$, which may be a value that was not used in the first step for data curation and hence training, and provides an initial guess to the NS $\pi_\gamma$. If the generative model well approximates the desired conditional distribution, $p(\cdot|\alpha)$, then the NS $\pi_\gamma$ will convergence rapidly and robustly. We now explain in detail the choice of generative ML used in this paper.

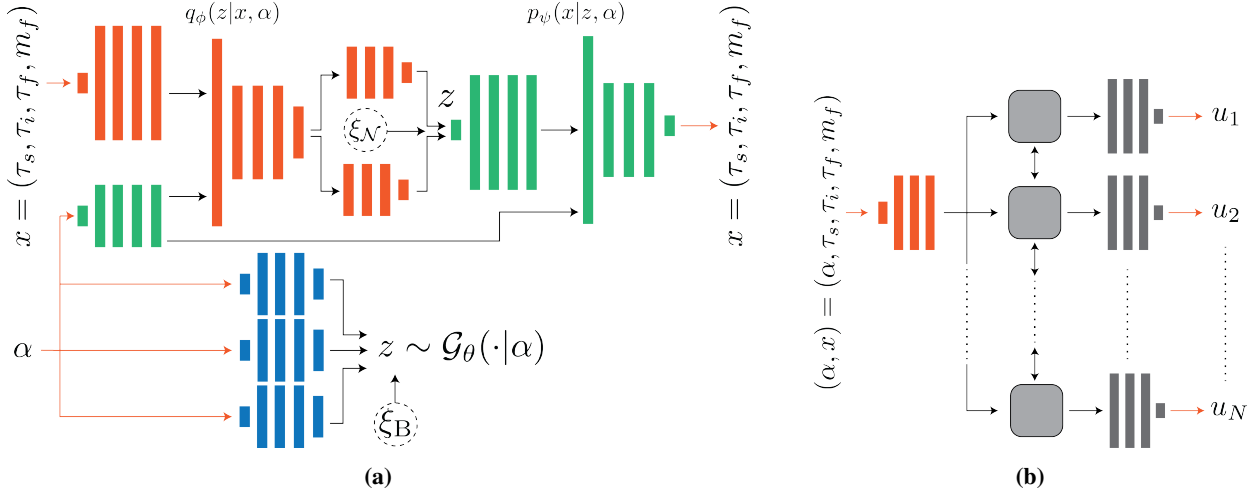

**Fig. 15** (a) The CVAE model with inputs and output corresponding to the nominal setup for Sec. V.B. $q_\phi(z|x,\alpha)$ is the approximate posterior (encoder), $\mathcal{G}_\theta(z|\alpha)$ the GMM prior, $z$ the latent variable, and $p_\psi(x|z,\alpha)$ the likelihood (decoder). (b) The LSTM model with inputs and output used in Sec. V.B, which generates the temporal control sequence $(u_1, \ldots, u_N)$.

### B. Conditional Variational Autoencoder and Gaussian Mixture Model

The generative ML model used in this work is a CVAE. Although the AmorGS framework shown in Fig. 14 proposes the ultimate goal of generating samples for the joint primal and dual space, $\mathcal{U}^h$, in the proceeding examples we restrict the CVAE to generating samples of just the primal space, $\mathcal{X}^h$. For initial guesses to the dual variables, the NS default (often a zero vector) is used (c.f., Sec. II.C).

A CVAE provides a tractable approach to generating samples $x \sim p(\cdot|\alpha)$ by first assuming that the data is partially observable. Hence there exists a latent variable $z$, such that the true underlying joint distribution is $p(x, z|\alpha)$. Using the standard definition of condition distributions, we then have $p(x, z|\alpha) = p(x|z, \alpha)p(z|\alpha)$. From this assumption, variational approaches may be introduced to learn a prior model, $p(z|\alpha)$, for the latent variable. The prior model is often assumed to be a distribution that is easy to sample from, such as a Gaussian distribution. To enable greater flexibility in learning a distribution with support on the $\alpha$-dependent funnel and clustering structures of our global search problems, we follow the clustering work of Jiang et al. [56] on images and Xiong et al. [57] on single-cell gene sequencing that provides the base research for variational autoencoders (VAE) and Gaussian Mixture Models (GMM). Similarly to Sohn et al. [58], we extend this work onto a conditional generative framework.

Figure 15 provides a detailed outline of the CVAE network architecture used in this work. The data input and samples generated, denote by $x$, are written in explicit form as to the data type used in the LT CR3BP problem. The CVAE consists of three components: the encoder $q_\phi(z|x, \alpha)$, shown as orange NN layers; the prior $\mathcal{G}_\theta(z|\alpha)$, shown as blue NN layers; and the likelihood $p_\psi(x|z, \alpha)$, shown as green NN layers. The CVAE network parameters that must be learned are denoted by $\phi, \theta$ and $\psi$. As implied from Step 2 to 3 of Fig. 14, although the encoder is needed for training, it is discarded for online usage of the CVAE.

We choose our prior to be a GMM with $K$-components, conditioned on the parameter $\alpha$,

$$p(z|\alpha) \equiv \mathcal{G}_\theta(z|\alpha) = \sum_{i=1}^{K} w_i(\alpha)\mathcal{N}_i(z; \boldsymbol{\mu}, \boldsymbol{\Sigma}|\alpha) = \mathcal{N}(z; \mu_i, \Sigma_i|i, \alpha)\mathrm{B}(i; w|\alpha), \tag{34}$$

where for each $\alpha$,

$$\sum_{i=1}^{K} w_i(\alpha) = 1, \qquad w_i(\alpha) \geq 0, \ \forall i \in \{1, \ldots, K\}. \tag{35}$$

In Eq. (34) we use $\mathcal{G}_\theta(z|\alpha)$ to clearly denote the conditional GMM, $z$ represents the latent variables, and $w_i(\alpha)$ are the weights of each of the Gaussian mixture components. For each Gaussian $\mathcal{N}_i$ in the mixture, $\boldsymbol{\mu}$ represents the mean and $\Sigma$ the covariance matrix. As is common in applications of CVAEs, we assume that the covariance matrix of each Gaussian component is a diagonal matrix. This simplifying assumption enables numerically efficient computation of determinants, which are required for the variational component of the loss function to be maximized. The final relation of Eq. (34) makes use of a generalized Bernoulli distribution over the indices $(1, \ldots, K)$, each index with weight $w_i$. The three NNs shown in blue in Fig. 15 generate the weights $w$, the means $\mu$, and the covariances $\Sigma$ of the GMM components.

The CVAE network architecture used in this paper makes use of multilayer perceptrons with a leaky rectified linear unit as the activation function in the hidden layers. A sigmoid function is used as the activation function before the final output $x$ and batch normalization is used during training. Tables 3, 4, and 5 provide a breakdown of the NN layer sizes used in the experimental problems. The Embed_x_layer, Encode_layer, Encode_$\mu$_layer and Encode_$\Sigma$_layer are shown as the encoder (in orange) in Fig. 15. The Embed_$\alpha$_layer, Embed_z_layer, and Decode_x_layer are shown as the decoder (in blue) in Fig. 15. Lastly, the GMM_w_layer, GMM_$\mu$_layer and GMM_$\Sigma$_layer are shown as the GMM (in green) in Fig. 15.

### C. Loss Function and Evidence Lower Bound

Assuming that our observed data $(x_i) \subseteq \mathcal{A}_{\alpha,\beta}$ are samples from the parameterized distribution $p_{\psi,\theta}(x|\alpha)$, the aim of training is to maximize the log-likelihood of $p_{\psi,\theta}(x|\alpha)$. By introducing the assumption of a latent variable and using the ideas of importance sampling and Jensen's inequality, it is easy to arrive at a lower bound on the log-likelihood that is computationally more tractable to maximize,

$$\log p_{\psi,\theta}(x|\alpha) = \log \mathbb{E}_{q_\phi(z|x,\alpha)} \left[ \frac{p_{\psi,\theta}(x,z|\alpha)}{q_\phi(z|x,\alpha)} \right] \geq \mathbb{E}_{q_\phi(z|x,\alpha)} \left[ \log p_\psi(x|z,\alpha) + \log p_\theta(z|\alpha) - \log q_\phi(z|x,\alpha) \right]$$

$$= \mathbb{E}_{q_\phi(z|x,\alpha)} \left[ \log p_\psi(x|z,\alpha) \right] - \mathrm{D}_{\mathrm{KL}} \left( q_\phi(z|x,\alpha) || p_\theta(z|\alpha) \right). \tag{36}$$

The lower bound in Eq. (36) is known as the Evidence Lower Bound (ELBO). The first term amounts to a reconstruction error, and assuming that the likelihood is a Gaussian distribution, is proportional to the Mean Squared Error (MSE) between the original and reconstructed data. The second term is the Kullback-Liebler (KL)-divergence between the posterior distribution $q_\phi(z|x,\alpha)$ and the prior $p_\theta(z|\alpha)$. The encoder layers of the CVAE are chosen to output the mean and covariance parameters for a Gaussian distribution; this is therefore the posterior distribution conditioned on the original data and parameter $\alpha$. The prior is chosen to be our GMM, $p_\theta(z|\alpha) = \mathcal{G}_\theta(z|\alpha)$. Hence, minimizing the KL-divergence promotes choosing parameters for the prior such that the statistical distance between these two distributions are close.

To enable backpropagation through the posterior and to enable computational tractability of computing the KL-divergence, the standard generative ML re-parameterization trick is used, whereby we introduce random samples from a standard Gaussian $\xi_N$ and generalized Bernoulli $\xi_B$ as shown in Figs. 14 and 15. The generalized Bernoulli random variable provides the selection of a single component of the GMM and hence an analytic solution to the KL-divergence. We refer the reader to Kingma and Welling [59] for details on the re-parameterization trick.

As is similarly done in Sohn et al. [58], we introduce a weighting parameter $\eta_\mathcal{L} > 0$ to re-balance the emphasis of the two components in ELBO, giving the final loss function to be maximized in the form,

$$\max_{\psi,\phi,\theta} \left\{ \mathcal{L}_{\text{CVAE}}(\psi,\phi,\theta) \equiv \mathbb{E}_{q_\phi(z|x,\alpha)}\left[\log p_\psi(x|z,\alpha)\right] - \eta_\mathcal{L}\, \text{D}_{\text{KL}}\left(q_\phi(z|x,\alpha)||\mathcal{G}_\theta(z|\alpha)\right) \right\}. \tag{37}$$

In the proceeding problems, we use single realization sampling to approximate the expectations in the loss function and set $\eta_\mathcal{L} = 1\text{E-}4$.

**D. Long Short-Term Memory**

The control variables $(u_1, \ldots, u_N)$ in the transcription of the LT CR3BP problem possess temporal correlation, which causes difficulty in capturing the structure of the solutions using naive approaches; this will be demonstrated from the numerical results of the ablation study given in Sec. V.B.7. Therefore in the LT CR3BP problem, the CVAE is used to generate samples of the time variables $(\tau_s, \tau_i, \tau_f)$ and final mass $m_f$, and an LSTM model is used to generate the controls $(u_1, \ldots, u_N)$ conditioned on the time variables, final mass, and $\alpha$. This idea is shown in Fig. 15.

An LSTM [60] is a special type of Recurrent Neural Network (RNN), which has a specific architecture that is well-suited for dealing with sequential data. An RNN usually has hidden layers that form directed cycles, which effectively enables it to store information about past inputs in its internal state, therefore allowing it to exploit temporal dependencies within the data sequence. A key feature of the LSTM architecture is its ability to alleviate the vanishing and exploding gradient problem encountered in traditional RNNs. LSTMs have the ability to selectively remember or forget information through a series of gating mechanisms, which makes them particularly effective for handling longer sequentially correlated data streams. The LSTM network architecture shown in Fig. 15 has an encoder layer, which feeds into the LSTM bidirectional components, each of which is connected to a three layer decoder (see Graves and Schmidhuber [61] for details on the bidirectional architecture). An MSE loss function is used for training the LSTM.

Table 3  The CVAE network architecture for the De Jong's 5th function problem.

| CVAE | | | |
| --- | --- | --- | --- |
| **Layer Name** | **Layer Size** | **Layer Name** | **Layer Size** |
| Embed_$x$_layer | [2, 32, 64, 64] | Embed_$\alpha$_layer | [1, 32, 64, 64] |
| Encode_layer | [128, 64, 64] | Encode_$\mu$_layer | [64, 32, 2] |
| Encode_$\Sigma$_layer | [64, 32, 2] | Embed_$z$_layer | [2, 32, 64, 64] |
| Decode_$x$_layer | [128, 64, 64, 2] | GMM_$w$_layer | [1, 512, 512, 512, 2] |
| GMM_$\Sigma$_layer | [1, 512, 512, 512, 4] | GMM_$\mu$_layer | [1, 512, 512, 512, 4] |

## V. Test Results

### A. De Jong's Fifth Function

In this section we present the results of applying the AmorGS framework to the De Jong's fifth function example, which was formulated in Sec. III.A, to build an intuition of the framework application before a deeper analysis of the main LT CR3BP problem. In particular, we aim to test the performance of predicting the positions and rotations of the minima clusters when an *a priori* unknown rotation $\alpha$ is considered. To mimic a clustering pattern (e.g., as seen in the LT CR3BP optimization problem), we set the matrix $A$ in Eq. (20) as,

$$A = \begin{pmatrix} -32 & -32 & -28 & -28 & 12 & 12 & 18 & 18 \\ 32 & 28 & 32 & 28 & -12 & -18 & -12 & -18 \end{pmatrix}, \tag{38}$$

such that the eight local minima are grouped into two clusters, each with four minima.

The data curation phase included the collection of 80,000 numerical local minima solutions solved to an optimality tolerance of 1E-6 with a Broyden-Fletcher-Goldfarb-Shanno (BFGS) method. These solutions were collected from 10,000 values of $\alpha$ sampled uniformly from 51 evenly spaced points in $[0, \pi/2]$ (i.e., $\alpha \sim U((i\pi/100)_{i=0}^{50})$). Using this data, a CVAE with the network architecture given in Table 3 was trained using PyTorch [62]. The GMM prior was chosen to consist of two components, each a 2-dimensional Gaussian with a diagonal covariance matrix.

Examples of the training data are shown in Fig. 16a for different values of $\alpha$, which shows the changing position of the local minima and rotation of the clusters as the $\alpha$-parameter varies. In Fig. 16b we show the prediction results of the trained CVAE in the AmorGS framework for values of $\alpha$ not used in the training data set. In particular, we sample the trained CVAE for 1,000 samples, which we show with blue markers. Since the analytic form of the minima is given by Eq. (20) with matrix $A$ from Eq. (38), we are able to superimpose the true minima, which are shown with orange

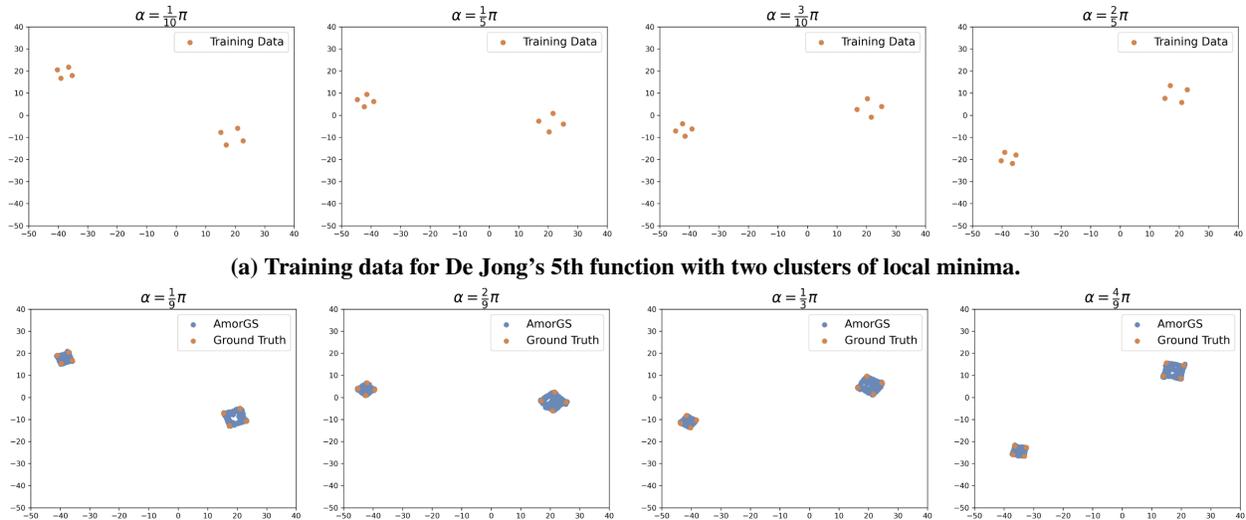

(a) Training data for De Jong's 5th function with two clusters of local minima.

(b) 1,000 samples from the AmorGS framework for $\alpha \in \{\frac{1}{9}\pi, \frac{2}{9}\pi, \frac{1}{3}\pi, \frac{4}{9}\pi\}$, which are all parameter values not in the training set for the De Jong's 5th function.

Fig. 16   Training and testing/prediction data for the De Jong's 5th function.

markers in Fig. 16b. A neighborhood of the clusters are clearly covered by the samples in each new case of $\alpha$, which indicates that the learned distribution largely has support over the local neighborhoods of the minima. There was a low probability of sampling outside of these local neighborhoods, which was typically only seen if the number of samples was increased by another order of magnitude.

## B. Low Thrust Circular Restricted Three-Body Problem

This section begins by providing details on the experimental setup of the LT CR3BP problem, details of the network architectures, before providing analysis on the ability of the learned models to predict the solution topology for parameter cases not in the training data set. As in Sec. III.B, this analysis first looks at the prediction of the time variables, then final mass and other control variables. An analysis of training on reduced data is explored. A thorough study of the performance on warm starting the NS using the trained framework is detailed, including an ablation study where various components of the framework are removed to demonstrate their individual impact on the overall performance. Lastly, we demonstrate the support of the learned distribution on the neighborhoods of the funnel and cluster structure such that the NS produces solutions that replicate the diversity of the full optimal solution set.

### 1. Experiment Setup

For the LT CR3BP problem, the $\alpha$-parameter represents a scaling of the maximum allowable thrust, as shown in Eq. (25). In the data curation step, 300,000 solutions are collected from the union of sets $\mathcal{A}_{\alpha,\beta}$ with $\alpha$ taken over 12 values,

$$\alpha \in \{i/10 \mid i \in \{1, 2, \ldots, 10\}\} \cup \{0.13, 0.16\}, \tag{39}$$

Table 4  The CVAE and LSTM network architecture for the LT CR3BP transfer problem.

| CVAE | | | |
|---|---|---|---|
| Layer Name | Layer Size | Layer Name | Layer Size |
| Embed_$x$_layer | [4, 1024, 1024, 1024, 1024] | Embed_$\alpha$_layer | [1, 256, 256, 256, 256] |
| Encode_layer | [1280, 512, 512, 512, 128] | Encode_$\mu$_layer | [128, 128, 128, 4] |
| Encode_$\Sigma$_layer | [128, 128, 128, 4] | Embed_$z$_layer | [4, 1024, 1024, 1024, 1024] |
| Decode_$x$_layer | [1280, 512, 512, 512, 4] | GMM_$w$_layer | [1, 512, 512, 512, 20] |
| GMM_$\mu$_layer | [1, 512, 512, 512, 80] | GMM_$\Sigma$_layer | [1, 512, 512, 512, 80] |
| LSTM | | | |
| Layer Name | Layer Size | Layer Name | Layer Size (For each time step) |
| Encoder_layer | [5, 512, 512, 512] | Decoder_layer | [512, 512, 512, 3] |

and with $\beta$ = 415 kg acting as the threshold designating high-quality solutions. Approximately 25,000 solutions are generated for each value of $\alpha$ given in Eq. (39), with slightly more for the lower values of $\alpha$, which represent the more difficult LT transfers. To generate these solutions, initial guesses consisting of samples from uniform distributions over $\tau_s \in [0, 40]$ TU, $\tau_i, \tau_f \in [0, 15]$ TU, and $m_f \in [350, 450]$ kg, as well as $\boldsymbol{u}_i \in [0, 2\pi] \times [-\pi/2, \pi/2] \times [0, T_{\max}]$. To thoroughly exploit the solution space, the maximum run time of the solver $\pi_\gamma$ was set to 500 seconds and the major iteration limit was set to 1,000. A feasibility and optimality tolerance of 1E-3 was used.

The Adam optimizer [63] was used for the training of the CVAE and the LSTM. The CVAE was trained for 1200 epochs and the LSTM for 600 epochs. These values were sufficiently high such that full convergence was obtained for the two models. The training process was completed on a single NVIDIA A100 GPU and required approximately 2 hours.

Validation of the AmorGS framework is done for two parameter values not in the set given by Eq. (39); in particular $\alpha \in \{0.15, 0.85\}$. The first of these values has data near it (in 0.13 and 0.16), whereas the latter has a greater spread in available data (closest data points being 0.8 and 0.9). The lower value of $\alpha = 0.15$ corresponds to a lower maximum available thrust, and corresponds to a regime where the solution topology changes more dramatically per equal step in $\alpha$ in comparison to the higher thrust case. Therefore using training data closer to the unseen 0.15 case is possibly better for training that case, though to what degree is difficult to quantify.

*2. Model Architecture*

The architecture of the CVAE and LSTM models are shown in Table 4. Every layer in Table 4 is a fully connected layer and the size indicates its width and length. The GMM prior is chosen to consist of $K = 20$ components, each of which is a 4-dimensional Gaussian with a diagonal covariance matrix. The models are built and trained concurrently with PyTorch [62].

Table 5    The vanilla CVAE network architecture to baseline the LT CR3BP transfer problem.

| Vanilla CVAE | | | |
|---|---|---|---|
| **Layer Name** | **Layer Size** | **Layer Name** | **Layer Size** |
| Embed_$x$_layer | [64, 1024, 1024, 1024, 1024] | Embed_$\alpha$_layer | [1, 256, 256, 256, 256] |
| Encode_layer | [1280, 512, 512, 512, 128] | Encode_$\mu$_layer | [128, 128, 128, 64] |
| Encode_$\Sigma$_layer | [128, 128, 128, 64] | Embed_$z$_layer | [64, 1024, 1024, 1024, 1024] |
| Decode_$x$_layer | [1280, 512, 512, 512, 64] | | |

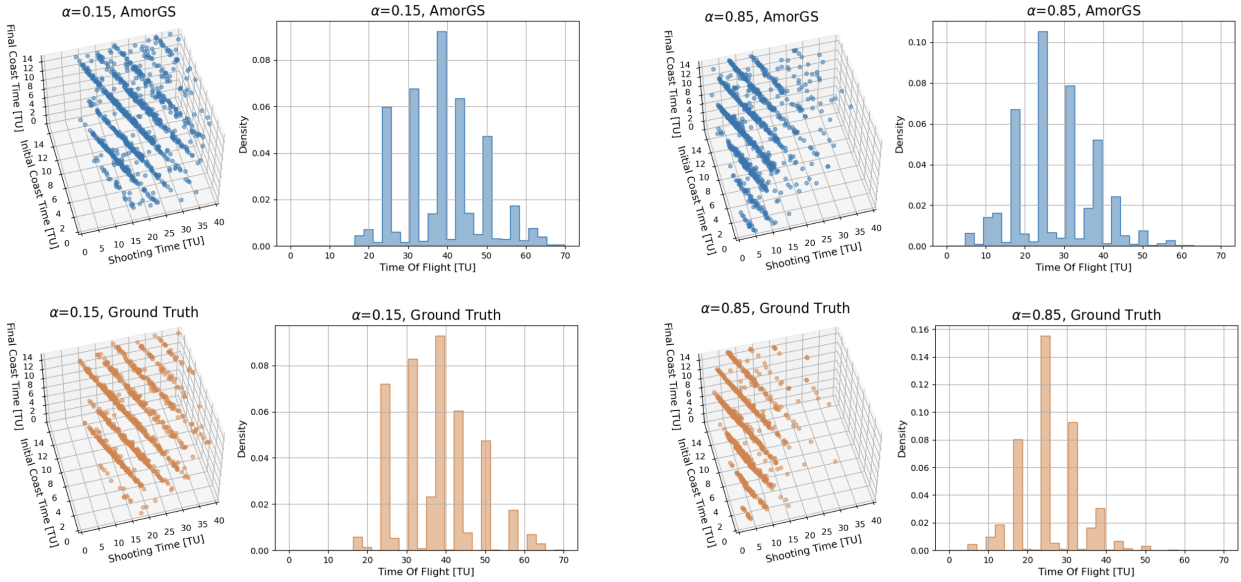

**Fig. 17    Comparison of the time variable prediction between the AmorGS framework and ground truth data. (Left) images display the $\alpha = 0.15$ case and (Right) the $\alpha = 0.85$ case.**

*3. Prediction of Time Variables*

Figure 17 demonstrates the ability of the learned AmorGS framework to predict the structure of the time variables $\tau$ in the $\alpha \in \{0.15, 0.85\}$ case. The top row of the figure shows the predicted (sampled) values from the learned framework (in blue), whereas the bottom row includes the ground truth data for comparison (in orange). Histograms of the total time-of-flight are included, which provides another comparison of the closeness of the learned conditional distribution for these time variables in comparison to the true structure. The modes of the distribution are largely well captured, and the translation of the hyperplane structure seen in the time variables is well predicted by the learned conditional distribution. In the $\alpha = 0.85$ case, some higher time-of-flight, but low probability modes are predicted, yet not seen in the true data. The slight performance improvement in the $\alpha = 0.15$ case versus the 0.85 case may be due to the training data being closer to the 0.15 scenario.

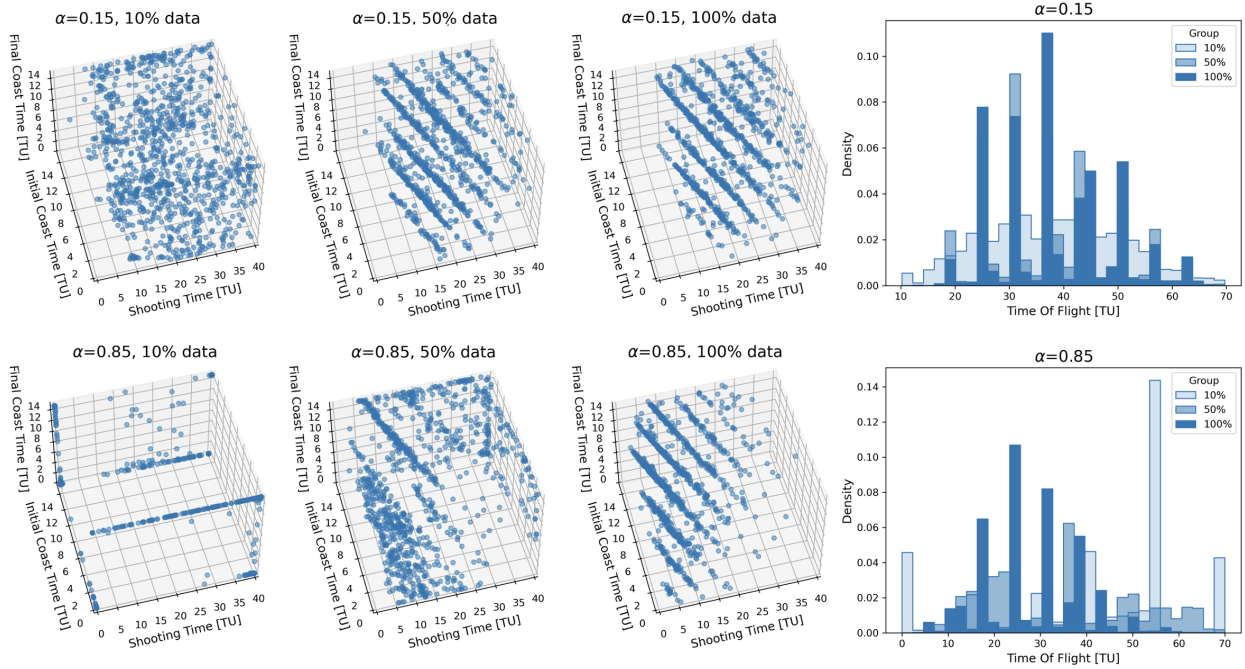

**Fig. 18   Hyperplane predictions of time variables $\tau$ for $\alpha \in \{0.15, 0.85\}$ using 10%, 50%, and 100% of the total training data.**

*4. Reduced Training Dataset*

Achieving good predictions of the solution topology based on a limited training data set is desirable, especially for high fidelity problems that may take substantial time to solve. In this section, we investigate the effect of reducing the training data set size on the prediction performance of the CVAE model for the time variables.

Figure 18 shows the hyperplane predictions of the time variables, $\tau$, using 10%, 50%, and 100% of the training data, as well as a subfigure providing comparisons of histograms for the reduced data cases. From Fig. 17 it is clear that the 100% data case was sufficient for the model to sample from the neighborhood of the true data hyperplane. In the case of $\alpha = 0.15$, the CVAE model is still capturing the hyperplane structure at the 50% data case, but now with a slight increase in diffusion in the direction normal to the hyperplanes. This is also seen in the histograms of Fig. 18. Reducing the amount of training data to the 10% case results in ineffective learning of the hyperplane structure, with a broadly diffuse sampling over the parameter ranges for the $\tau$ variables. The behavior in the $\alpha = 0.85$ case is more dramatic, and already at 50% of training data most hyperplanes are no longer captured accurately. The exception is a hyperplane corresponding to time-of-flight near 37.5 TU. A collapse of the condition distribution in the 10% data case is qualitatively different than the $\alpha = 0.15$ case.

Because all $\alpha$ scenarios from Eq. (39) are used in the training phase, it is hypothesized that the contribution of the $\alpha = 0.85$ scenario is diminished in the approximation of the loss function given in Eq. (37) and its gradient. The behavior of this collapse was replicate for independent training experiments and may indicate that the extrema reached for the loss function is broad and stable for the given network architecture. One possible approach to alleviate the collapse in the case of abundant data for select values of $\alpha$, but fewer for new trials, is to consider the base learned model (e.g., the 100% data case considered here) as a foundational model, and leverage techniques to fine-tune the model for the reduced data scenarios. Future efforts will further investigate the collapse phenomena and applicability of the foundational model paradigm.

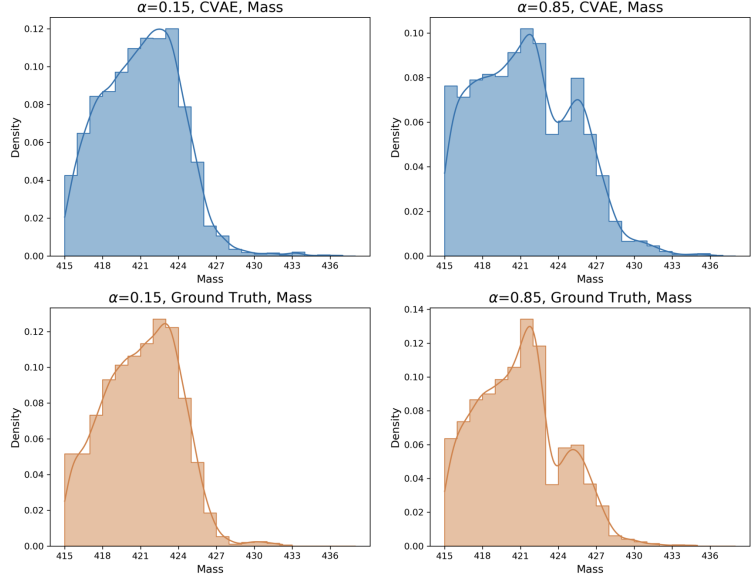

**Fig. 19** **Histograms of the mass variable $m_f$ for $\alpha \in \{0.15, 0.85\}$ between CVAE predictions and the ground truth.**

*5. Prediction of the Mass Variable*

Figure 19 provides a histogram of 4,000 samples of the mass variable $m_f$ from the learned CVAE model with comparison to samples of the ground truth data in the cases of $\alpha \in \{0.15, 0.85\}$. A kernel density estimation is overlayed on each histogram plot. Both cases of $\alpha$ are well approximated by the CVAE model, though the $\alpha = 0.15$ case is again superior. The case of $\alpha = 0.85$ predicts a greater number of solutions with lower delivered final mass than is seen in the ground truth data. A nice feature that is duplicated in the $\alpha = 0.85$ case is the second mode appearing around 425 kg, whereas the $\alpha = 0.15$ is very much a unimodal marginal distribution for final mass.

*6. Prediction of the Control Variables*

To analyze the prediction capabilities of the control variables $(\boldsymbol{u}_1, \ldots, \boldsymbol{u}_N)$, the samples from the LSTM and ground truth data for the $\alpha \in \{0.15, 0.85\}$ cases are converted from spherical coordinates to cartesian. Figure 20 shows the mean and 95% confidence interval of the control components $u_x$, $u_y$, and $u_z$ across the 20 control segments based on 4,000 samples from the LSTM, ground truth data, and a sampling from a uniform distribution. The sampling from the uniform distribution is simple to provide a baseline to judge the variability for the 4,000 samples. As can be expected for this problem, which entails multiple revolutions about Earth as the spacecraft transits from the initial boundary condition at the end of a GTO LT spiral to the terminal boundary condition on the invariant manifold to the $L_1$ Halo, the control components in the $u_x, u_y$ directions of the CR3BP canonical frame largely average out. The trajectories do need to perform out-of-plane changes and hence the distinctive presence of an on average positive $u_z$ component. It is likely that due to the closeness of the training data for $\alpha = 0.16$ to the test case of $\alpha = 0.15$, the $\alpha = 0.15$ predictions are able to better follow the average trend in $u_z$ as seen in Fig. 20. As evidenced by the effectiveness in warm starting the NS $\pi_\gamma$, to be shown in Sec. V.B.7, the prediction of the $u_z$ component for the $\alpha = 0.85$ case still appears largely satisfactory.

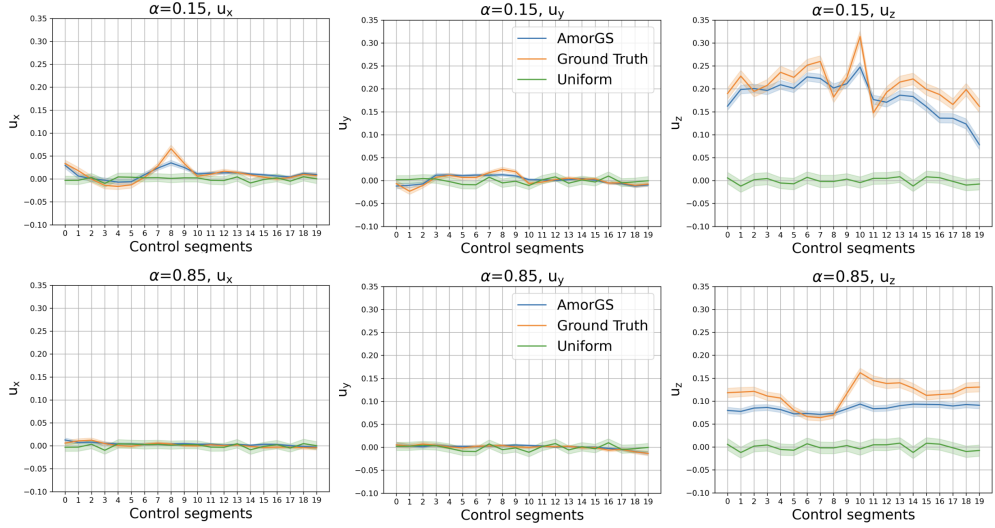

**Fig. 20** **95% confidence interval of control variables in cartesian coordinates for $\alpha \in \{0.15, 0.85\}$ among LSTM predictions, ground truth, and samples from a uniform distribution.**

*7. Warm Start Performance and Ablation Study*

As detailed in the preceding sections, the AmorGS framework provides samples in the neighborhood of the desired minima set $\mathcal{A}_{\alpha,\beta}$ for the cases of $\alpha \in \{0.15, 0.85\}$, which is out of the trained distribution from Eq. (39). In this section we quantify the ability for these predictions to improve the rate of successful convergence when solved using the NS, and the ability to reduce the computational time in doing so. To benchmark the factor of improvement in each of these categories, and to understand the role of each component in the AmorGS framework, we perform an ablation study. In particular, we consider five methods of sampling:

- **Uniform**: Time $\tau$, final mass $m_f$, and control $(\boldsymbol{u}_1, \ldots, \boldsymbol{u}_N)$ variables are all sampled from uniform distributions on their respective domains.

- **CVAE time & mass, uniform control**: Time $\tau$ and final mass $m_f$ variables are sampled from the CVAE model. Control $(\boldsymbol{u}_1, \ldots, \boldsymbol{u}_N)$ variables are sampled from a uniform distribution.

- **Uniform time & mass, LSTM control**: Time $\tau$ and final mass $m_f$ variables are sampled from uniform distributions. Control $(u_1, \ldots, u_N)$ variables are sampled from the LSTM model.
- **Vanilla CVAE**: Time $\tau$, final mass $m_f$, and control $(u_1, \ldots, u_N)$ variables are all jointly sampled from a vanilla CVAE that uses a 64-dimensional Gaussian as the latent distribution.
- **AmorGS**: Time $\tau$, final mass $m_f$, and control $(u_1, \ldots, u_N)$ variables are all sampled from the full AmorGS framework (i.e., CVAE and LSTM models).

The **Uniform** case is the naive approach that was used to generate the training data and would presumably be used for the cases $\alpha \in \{0.15, 0.85\}$ given no prior knowledge. The **CVAE time & mass, uniform control** case provides an understanding of the effectiveness of using the LSTM for the temporally correlated control variables $(u_1, \ldots, u_N)$. The **Vanilla CVAE** will further clarify the effectiveness of the LSTM by showing that the temporal correlation structure does indeed warrant special attention. Lastly, the **Uniform time & mass, LSTM control** case provides similar information regarding the ability for the CVAE with GMM prior to capture the structure of the time and final mass variables. All versions of CVAE and LSTM are trained on the same data sets.

Table 6 provides the summarized results of testing each of the aforementioned approaches. In particular, 200 samples were drawn for each approach and used as initial guesses to the NS $\pi_\gamma$. As shown by the first row in each of the tables for the $\alpha = 0.15$ and $0.85$ cases, the percentage of initial guesses that converged within the allowable 500 seconds or 1,000 major iterations of $\pi_\gamma$ was significantly higher for the full AmorGS framework; especially in the $\alpha = 0.15$ case, which is more computationally intensive and sensitive to initial conditions due to the longer time-of-flight. The reduction in percentage of converged solutions from the $\alpha = 0.85$ to $0.15$ case is largest in the methods that did not use an LSTM. To be clear and explicit, we see a drop in the percentage of converged solutions for the **Uniform**, **CVAE time & mass, uniform control**, and **Vanilla CVAE** methods of 14%, 15.5%, and 24.5% respectively. But the **Uniform time & mass, LSTM control** and **AmorGS** methods only drop by 11% and 11.5% respectively. Yet, to gain the full power of the LSTM, a good sample of the time and final mass parameters are required, hence the improvement of the **AmorGS** method versus the **Uniform time & mass, LSTM control**.

Table 6 also provides the solver solution time for the converged cases. The **AmorGS** method is superior in each of the case studies and in comparison to any other approach. In comparison to the **Uniform** method, it is on average (mean) between 1.5 to 2.5 times quicker. The comparison of the speed improvement in the median solution is 2.5 to 6 times quicker.

Figure 21 provides histograms for required solver time to convergence for the solutions given in Table 6, and for each the five methods. The final column of the histograms, corresponding to the full AmorGS framework, shows a clear mode at the lower solver times. Additionally, as is also reported in Table 6, only the AmorGS framework had samples that the solver was able to converge in less than 10 seconds. The AmorGS framework was able to do this for 4% of the converged samples.

Table 6  First row of each table: the percentage of 200 initializations that converge to the solutions in $\mathcal{A}_{\alpha,\beta}$ for $\alpha \in \{0.15, 0.85\}$ and $\beta = 415$ kg. Remaining rows: solver time statistics for converged solutions (mean, minimum, 25% percentile, and 50% (median) percentile).

(a) $\alpha = 0.15$

|  | Uniform | Uniform time & mass, LSTM control | CVAE time & mass, Uniform control | Vanilla CVAE | AmorGS |
|---|---|---|---|---|---|
| Percentage (%) | 28 | 29.5 | 42 | 30 | **62.5** |
| $\pi_\gamma$ solve time (s) | | | | | |
| Mean | 195.48 | 196.17 | 178.90 | 201.78 | **112.03** |
| Minimum | 25.08 | 12.84 | 20.06 | 32.80 | **4.37** |
| 25% | 103.29 | 89.89 | 79.24 | 97.03 | **25.89** |
| 50% (Median) | 169.31 | 168.56 | 152.27 | 204.38 | **64.14** |

(b) $\alpha = 0.85$

|  | Uniform | Uniform time & mass, LSTM control | CVAE time & mass, Uniform control | Vanilla CVAE | AmorGS |
|---|---|---|---|---|---|
| Percentage (%) | 42 | 40.5 | 57.5 | 54.5 | **74** |
| $\pi_\gamma$ solve time (s) | | | | | |
| Mean | 152.28 | 152.96 | 113.73 | 143.17 | **68.38** |
| Minimum | 10.16 | 10.35 | 11.19 | 12.20 | **3.02** |
| 25% | 52.23 | 66.16 | 30.96 | 54.81 | **12.56** |
| 50% (Median) | 121.48 | 128.38 | 75.26 | 110.75 | **24.32** |

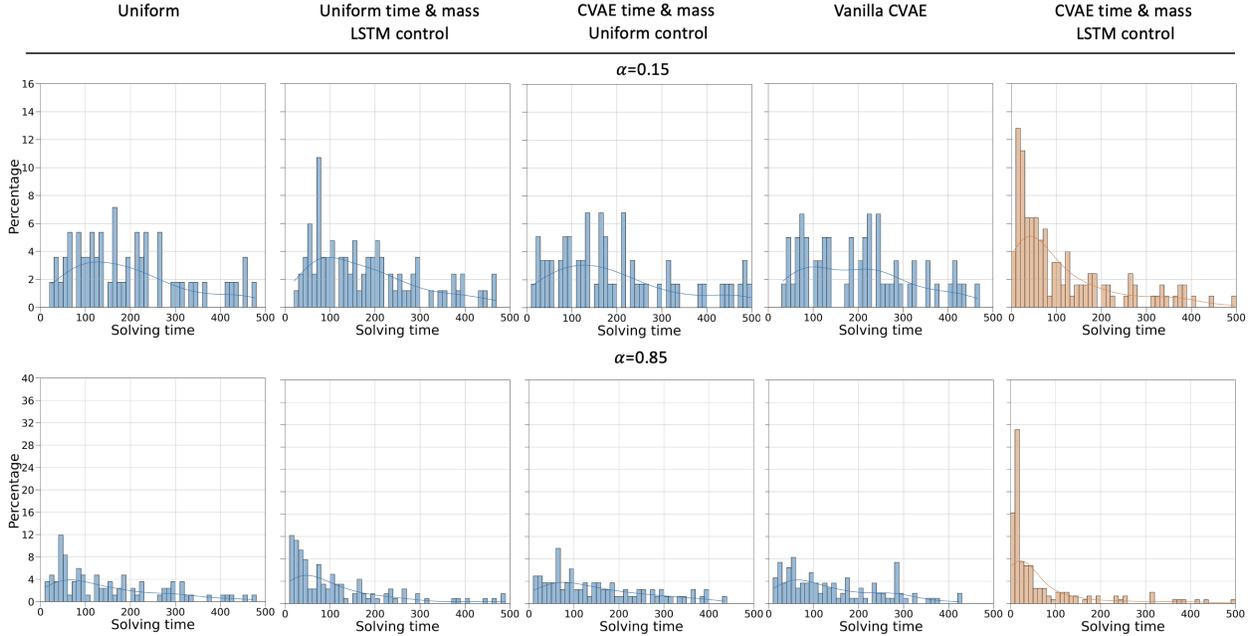

Fig. 21  Histograms of solution time for the NS $\pi_\gamma$ for the cases of $\alpha \in \{0.15, 0.85\}$. 200 samples from the five methods described in Sec. V.B.7. A kernel density estimate is overlayed for each histogram.

Since the global search problem requires a practitioner to trade time spent on exploration versus local exploitation (i.e., solving with $\pi_\gamma$), it is common to set the run time of the NS at a value that guarantees most good initial guesses will converge, but excessive NS run time is avoided[§]. The median time to converge for the AmorGS method is therefore taken as an upper limit run time in the next study. As reported in Table 6, we take the maximum run time for the $\alpha = 0.15$ case to be 64 seconds, and in the $\alpha = 0.85$ case to be 24 seconds. Figure 22 then shows the percentage of 4,000 samples taken for each method that converge within these maximum time limits. Note that the median reported in Table 6 is for converged solutions and therefore the success rate of 33 and 35% for the $\alpha = 0.15$ and 0.85 cases respectively for the AmorGS method is expected. This also confirms that the sample size of 200 in Table 6 was sufficient. Figure 22 makes clear that with an appropriately set maximum run time for the NS, a global search for the unseen cases of $\alpha \in \{0.15, 0.85\}$ can now be carried out at a rate of about ten times quicker than in the naive uniform search approach. This figure is also instructive in revealing the fact that identifying the solution structure with the CVAE is an important first step, but ultimately must be paired with the LSTM model to achieve dramatic search acceleration.

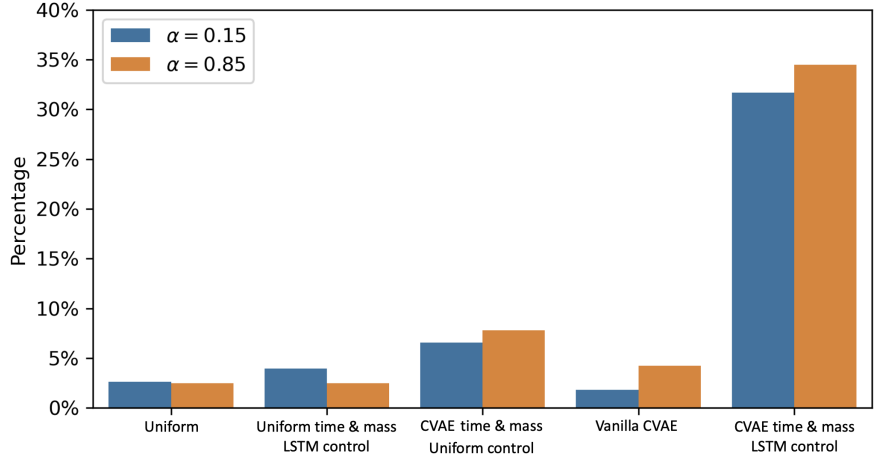

Fig. 22  **Percentage of 4,000 samples from each method that lead to numerically converged solutions in $\mathcal{A}_\alpha$ when the NS is limited to at most 64 seconds (in the $\alpha = 0.15$ case) and 24 seconds (in the $\alpha = 0.85$ case).**

*8. Diversity of Solutions*

A final analysis is provide for the AmorGS framework on the LT CR3BP test problem to confirm that the NS converges samples that are still diverse for each of the $\alpha = 0.15$ and 0.85 cases. This is important to consider, since an accelerated solution time, as detailed in the prior sections, is only meaningful if the converged solutions are diverse and empirically cover the set $\mathcal{A}_{\alpha,\beta}$ of interest. Figure 23 shows that this is indeed the case. The left column of this figure shows samples drawn from the learned CVAE and LSTM models for the AmorGS framework, which are then solved by the NS in the middle column with a limit on the solver run time. The right column of this figure then shows the ground truth data from a more extensive search of the problem for comparison. Because of the limit on the run time of the NS, it is able to easily converge samples that are already close to the hyperplane structure, but less so for samples that are farther from specific hyperplanes; such as that seen for the longer time-of-flight cases in the $\alpha = 0.85$ scenario.

---

[§]numerical solvers, especially when supplied with analytic gradients, can have the tendency to slowly reach convergence of the feasibility and optimality tolerance

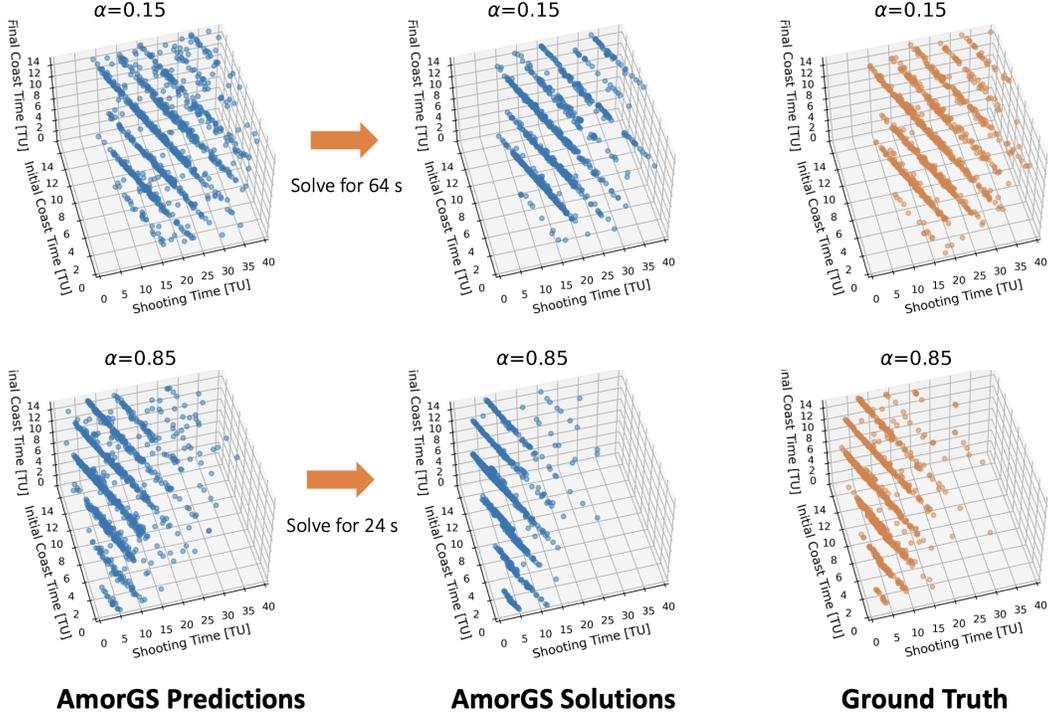

Fig. 23  Comparison of AmorGS samples for the unseen $\alpha \in \{0.15, 0.85\}$ **cases (left column), solution of these samples using a NS $\pi_\gamma$ with maximum run time of 64 and 24 seconds respectively (middle column), and the ground truth data for $\mathcal{A}_{\alpha,\beta}$ with $\beta = 415$ kg (right column).**

## VI. Conclusion

As the cadence of space mission design increases, and as high fidelity and more complex optimization is folded into earlier phases of the mission design process, it will become increasingly necessary to conduct more efficient global searches for optimal solutions. A similar statement could be made for the case of onboard autonomous space flight path planning. After formulating the global search problem for optimal trajectories as a problem in learning a conditional probability distribution, this paper demonstrated the effectiveness of using generative machine learning models within an amortized framework to provide high quality and diverse solution samples that are quickly converged by a numerical solver. These results are important in providing a first effective methodology to accelerating the long time duration low-thrust nonlinear optimal spacecraft trajectory problem using generative machine learning. A number of challenges still exist, but are placed in better focus from the mathematical framework and results of this paper. Of immediate interest is the extension to higher dimensional and multiple parameter cases, learning of $k$-neighborhoods, jointly learning primal and dual solution sets, and attempts to reduce the necessary data for effectively learning the conditional probability distribution for a given problem. Application of the framework to related low-thrust trajectory problems and study of the solution set structure using the techniques of Sec. III would also be informative for the wider community.

## Acknowledgments

The simulations presented in this article were performed on computational resources managed and supported by Princeton Research Computing, a consortium of groups including the Princeton Institute for Computational Science and Engineering (PICSciE) and the Office of Information Technology's High-Performance Computing Center and Visualization Laboratory at Princeton University.

## References

[1] Wales, D. J., and Doye, J. P. K., "Global Optimization by Basin-Hopping and the Lowest Energy Structures of Lennard-Jones Clusters Containing up to 110 Atoms," *The Journal of Physical Chemistry A*, Vol. 101, No. 28, 1997, pp. 5111–5116. https://doi.org/10.1021/jp970984n, URL https://doi.org/10.1021/jp970984n.

[2] Leary, R. H., "Global Optimization on Funneling Landscapes," *Journal of Global Optimization*, Vol. 18, 2000, pp. 367–383. https://doi.org/10.1023/A:1026500301312.

[3] Wales, D. J., "Exploring Energy Landscapes," *Annual Review of Physical Chemistry*, Vol. 69, No. 1, 2018, pp. 401–425. https://doi.org/10.1146/annurev-physchem-050317-021219, URL https://doi.org/10.1146/annurev-physchem-050317-021219, pMID: 29677468.

[4] Locatelli, M., and Schoen, F., *Global Optimization*, Society for Industrial and Applied Mathematics, Philadelphia, PA, 2013. https://doi.org/10.1137/1.9781611972672, URL https://epubs.siam.org/doi/abs/10.1137/1.9781611972672.

[5] Hartmann, J. W., Coverstone-Carroll, V. L., and Williams, S. N., "Optimal Interplanetary Spacecraft Trajectories via a Pareto Genetic Algorithm," *The Journal of the Astronautical Sciences*, Vol. 46, No. 3, 1998, pp. 267–282. https://doi.org/10.1007/BF03546237, URL https://doi.org/10.1007/BF03546237.

[6] Rauwolf, G. A., and Coverstone-Carroll, V. L., "Near-optimal low-thrust orbit transfers generated by a genetic algorithm," *Journal of Spacecraft and Rockets*, Vol. 33, No. 6, 1996, pp. 859–862. https://doi.org/10.2514/3.26850, URL https://doi.org/10.2514/3.26850.

[7] Russell, R. P., "Primer Vector Theory Applied to Global Low-Thrust Trade Studies," *Journal of Guidance, Control, and Dynamics*, Vol. 30, No. 2, 2007, pp. 460–472. https://doi.org/10.2514/1.22984, URL https://doi.org/10.2514/1.22984.

[8] Lawden, D., *Optimal Trajectories for Space Navigation*, Butterworths, London, 1963.

[9] Oshima, K., Campagnola, S., and Yanao, T., "Global search for low-thrust transfers to the Moon in the planar circular restricted three-body problem," *Celestial Mechanics and Dynamical Astronomy*, Vol. 128, No. 2, 2017, pp. 303–322. https://doi.org/10.1007/s10569-016-9748-2, URL https://doi.org/10.1007/s10569-016-9748-2.

[10] Lizia, P. D., and Radice, G., "Advanced Global Optimisation Tools for Mission Analysis and Design," FInal Report AO4532/18139/04/NL/MV, European Space Agency, 11 2004. URL https://www.esa.int/gsp/ACT/doc/ARI/ARI%20Study%20Report/ACT-RPT-MAD-ARI-03-4101b-GlobalOptimisation-Glasgow.pdf.

[11] Myatt, D. R., Becerra, V. M., Nasuto, S. J., and Bishop, J. M., "Advanced Global Optimisation Tools for Mission Analysis and Design," Final Report AO4532/18138.04/04/NL/MV, European Space Agency, 11 2004. URL https://www.esa.int/gsp/ACT/doc/ARI/ARI%20Study%20Report/ACT-RPT-MAD-ARI-03-4101a-GlobalOptimisation-Reading.pdf.

[12] Vasile, M., and Pascale, P. D., "Preliminary Design of Multiple Gravity-Assist Trajectories," *Journal of Spacecraft and Rockets*, Vol. 43, No. 4, 2006, pp. 794–805. https://doi.org/10.2514/1.17413, URL https://doi.org/10.2514/1.17413.

[13] Vasile, M., Minisci, E., and Locatelli, M., "Analysis of Some Global Optimization Algorithms for Space Trajectory Design," *Journal of Spacecraft and Rockets*, Vol. 47, No. 2, 2010, pp. 334–344. https://doi.org/10.2514/1.45742, URL https://doi.org/10.2514/1.45742.

[14] Reeves, C. R., and Yamada, T., "Genetic algorithms, path relinking, and the flowshop sequencing problem," *Evol. Comput.*, Vol. 6, No. 1, 1998, pp. 45–60. https://doi.org/10.1162/evco.1998.6.1.45, URL https://doi.org/10.1162/evco.1998.6.1.45.

[15] Addis, B., Cassioli, A., Locatelli, M., and Schoen, F., "A global optimization method for the design of space trajectories," *Computational Optimization and Applications*, Vol. 48, No. 3, 2011, pp. 635–652. https://doi.org/10.1007/s10589-009-9261-6, URL https://doi.org/10.1007/s10589-009-9261-6.

[16] European Space Agency Advanced Concepts Team, "Global Trajectory Optimisation Problems Database," , 2024. URL https://www.esa.int/gsp/ACT/projects/gtop/.

[17] Englander, J. A., and Englander, A. C., "Tuning Monotonic Basin Hopping: Improving the Efficiency of Stochastic Search as Applied to Low-Thrust Trajectory Optimization," *24th International Symposium on Space Flight Dynamics*, Laurel, Maryland, 2014. URL https://www.issfd.org/ISSFD_2014/ISSFD24_Paper_S7-3_Englander.pdf.

[18] Sims, J. A., and Flanagan, S. N., "Preliminary Design of Low-Thrust Interplanetary Missions," *AAS/AIAA Astrodynamics Specialist Conference*, Girdwood, Alaska, 1999.

[19] Englander, A. C., and Englander, J. A., "Walking the Filament of Feasibility," *AIAA/AAS Astrodynamics Specialist Meeting*, 2017. URL https://ntrs.nasa.gov/citations/20170008012.

[20] Englander, A., Englander, J., and Carter, M., "Hopping with an Adaptive Hop Probability Distribution," *AAS/AIAA Astrodynamics Specialist Conference*, South Lake Tahoe, CA, 2020.

[21] Ampatzis, C., and Izzo, D., "Machine Learning Techniques for Approximation of Objective Functions in Trajectory Optimisation," *IJCAI-09 Workshop on Artificial Intelligence*, 2009.

[22] Cassioli, A., Di Lorenzo, D., Locatelli, M., Schoen, F., and Sciandrone, M., "Machine learning for global optimization," *Computational Optimization and Applications*, Vol. 51, No. 1, 2012, pp. 279–303. https://doi.org/10.1007/s10589-010-9330-x, URL https://doi.org/10.1007/s10589-010-9330-x.

[23] Zhu, Y.-h., and Luo, Y.-Z., "Fast Evaluation of Low-Thrust Transfers via Multilayer Perceptions," *Journal of Guidance, Control, and Dynamics*, Vol. 42, No. 12, 2019, pp. 2627–2637. https://doi.org/10.2514/1.G004080, URL https://doi.org/10.2514/1.G004080.

[24] Stephens, C. P., and Baritompa, W., "Global Optimization Requires Global Information," *Journal of Optimization Theory and Applications*, Vol. 96, No. 3, 1998, pp. 575–588. https://doi.org/https://doi.org/10.1023/A:1022612511618, URL https://login.ezproxy.princeton.edu/login?url=https://www.proquest.com/scholarly-journals/global-optimization-requires-information/docview/196616878/se-2?accountid=13314.

[25] Schoen, F., *Two-Phase Methods for Global Optimization*, Springer US, Boston, MA, 2002, pp. 151–177. https://doi.org/10.1007/978-1-4757-5362-2_5, URL https://doi.org/10.1007/978-1-4757-5362-2_5.

[26] Bryson, A. E., and Ho, Y.-C., *Applied Optimal Control: Optimization, Estimation, and Control*, Hemisphere Publication Corporation, 1975.

[27] Liberzon, D., *Calculus of Variations and Optimal Control Theory: A Concise Introduction*, Princeton University Press, 2012. https://doi.org/10.2307/j.ctvcm4g0s, URL https://press.princeton.edu/books/hardcover/9780691151878/calculus-of-variations-and-optimal-control-theory.

[28] Prussing, J. E., *Primer Vector Theory and Applications*, Cambridge Aerospace Series, Cambridge University Press, 2010, pp. 16–36. https://doi.org/10.1017/CBO9780511778025.

[29] Hargraves, C., and Paris, S. W., "Direct Trajectory Optimization Using Nonlinear Programming and Collocation," *Journal of Guidance, Control, and Dynamics*, Vol. 10, No. 4, 1987, pp. 338–342. https://doi.org/10.2514/3.20223.

[30] Enright, P. J., and Conway, B. A., "Discrete Approximation to Optimal Trajectories Using Direct Transcription and Nonlinear Programming," *Journal of Guidance, Control, and Dynamics*, Vol. 15, No. 4, 1992, pp. 994–1002. https://doi.org/10.2514/3.20934.

[31] Betts, J. T., "Survey of Numerical Methods for Trajectory Optimization," *Journal of Guidance, Control, and Dynamics*, Vol. 21, No. 2, 1998, pp. 193–207. https://doi.org/10.2514/2.4231, URL https://doi.org/10.2514/2.4231.

[32] Ross, M., and Fahroo, F., "Legendre pseudospectral approximations of optimal control problems,", 2003. https://doi.org/10.1007/978-3-540-45056-6_21, URL http://hdl.handle.net/10945/66382.

[33] Betts, J. T., *Practical Methods for Optimal Control and Estimation Using Nonlinear Programming, Second Edition*, 2nd ed., Society for Industrial and Applied Mathematics, 2010. https://doi.org/10.1137/1.9780898718577, URL https://epubs.siam.org/doi/abs/10.1137/1.9780898718577.

[34] Conway, B. A., and Paris, S. W., *Spacecraft Trajectory Optimization Using Direct Transcription and Nonlinear Programming*, Cambridge Aerospace Series, Cambridge University Press, 2010, pp. 37–78. https://doi.org/10.1017/CBO9780511778025.

[35] Karush, W., "Minima of Functions of Several Variables with Inequalities as Side Constraints," Master's thesis, University of Chicago, Chicago, Illinois, 1939.

[36] Peressini, A. L., Sullivan, F. E., and J. J. Uhl, J., *The Mathematics of Nonlinear Programming*, Undergraduate Texts in Mathematics, Springer, 1988.

[37] Becker, R., and Lago, G., "A global optimization algorithm," *Proceedings of the 8th Allerton Conference on Circuits and Systems Theory*, Monticello, Illinois, 1970, pp. 3–12.

[38] Torn, A. A., "Cluster Analysis Using Seed Points and Density-Determined Hyperspheres as an Aid to Global Optimization," *IEEE Transactions on Systems, Man, and Cybernetics*, Vol. 7, No. 8, 1977, pp. 610–616. https://doi.org/10.1109/TSMC.1977.4309787.

[39] Rinnooy Kan, A. H. G., and Timmer, G. T., "Stochastic global optimization methods part I: Clustering methods," *Mathematical Programming*, Vol. 39, No. 1, 1987, pp. 27–56. https://doi.org/10.1007/BF02592070, URL https://doi.org/10.1007/BF02592070.

[40] Rinnooy Kan, A. H. G., and Timmer, G. T., "Stochastic global optimization methods part II: Multi level methods," *Mathematical Programming*, Vol. 39, No. 1, 1987, pp. 57–78. https://doi.org/10.1007/BF02592071, URL https://doi.org/10.1007/BF02592071.

[41] Locatelli, M., and Schoen, F., "Simple linkage: Analysis of a threshold-accepting global optimization method," *Journal of Global Optimization*, Vol. 9, No. 1, 1996, pp. 95–111. https://doi.org/10.1007/BF00121752, URL https://doi.org/10.1007/BF00121752.

[42] Locatelli, M., and Schoen, F., "Random Linkage: a family of acceptance/rejection algorithms for global optimisation," *Mathematical Programming*, Vol. 85, No. 2, 1999, pp. 379–396. https://doi.org/10.1007/s101070050062, URL https://doi.org/10.1007/s101070050062.

[43] Molga, M., and Smutnicki, C., "Test functions for optimization needs," , 2005. URL https://marksmannet.com/RobertMarks/Classes/ENGR5358/Papers/functions.pdf, access Oct. 18, 2024.

[44] Shekel, J., "Test Functions for Multimodal Search Techniques," *Fifth Annual Princeton Conference on Information Science and Systems*, 1971.

[45] Beeson, R., Sinha, A., Jagannatha, B., Bunce, D., and Carroll, D., "Dynamically Leveraged Automated Multibody (N) Trajectory Optimization," *AAS/AIAA Space Flight Mechanics Conference*, American Astronautical Society, Charlotte, NC, 2022.

[46] Gill, P. E., Murray, W., and Saunders, M. A., "SNOPT: An SQP Algorithm for Large-Scale Constrained Optimization," *SIAM Review*, Vol. 47, No. 1, 2005, pp. 99–131. https://doi.org/10.1137/S0036144504446096, URL https://doi.org/10.1137/S0036144504446096.

[47] Hunter, J. D., "Matplotlib: A 2D graphics environment," *Computing in Science & Engineering*, Vol. 9, No. 3, 2007, pp. 90–95. https://doi.org/10.1109/MCSE.2007.55.

[48] Li, A., Sinha, A., and Beeson, R., "Amortized Global Search for Efficient Preliminary Trajectory Design with Deep Generative Models," *AAS/AIAA Astrodynamics Specialist Conference*, Big Sky, Montana, 2023. URL http://arxiv.org/abs/2308.03960.

[49] Goodfellow, I. J., Pouget-Abadie, J., Mirza, M., Xu, B., Warde-Farley, D., Ozair, S., Courville, A., and Bengio, Y., "Generative Adversarial Nets," *Proceedings of the International Conference on Neural Information Processing Systems*, 2014, pp. 2672–2680. https://doi.org/10.48550/arXiv.1406.2661, URL https://proceedings.neurips.cc/paper_files/paper/2014/file/5ca3e9b122f61f8f06494c97b1afccf3-Paper.pdf.

[50] Sohl-Dickstein, J., Weiss, E., Maheswaranathan, N., and Ganguli, S., "Deep Unsupervised Learning using Nonequilibrium Thermodynamics," *Proceedings of the 32nd International Conference on Machine Learning*, Vol. 37, 2015, pp. 2256–2265. https://doi.org/10.48550/arXiv.1503.03585, URL http://proceedings.mlr.press/v37/sohl-dickstein15.pdf.

[51] Song, Y., Sohl-Dickstein, J., Kingma, D. P., Kumar, A., Ermon, S., and Poole, B., "Score-Based Generative Modeling through Stochastic Differential Equations," , 2020. https://doi.org/10.48550/arXiv.2011.13456, URL https://arxiv.org/abs/2011.13456.

[52] Ho, J., Jain, A., and Abbeel, P., "Denoising Diffusion Probabilistic Models," *34th Conference on Neural Information Processing Systems (NeurIPS 2020)*, Vol. 34, Vancouver, Canada, 2020. https://doi.org/10.48550/arXiv.2006.11239, URL https://proceedings.neurips.cc/paper/2020/file/4c5bcfec8584af0d967f1ab10179ca4b-Paper.pdf.

[53] Gräbner, J., Li, A., Sinha, A., and Beeson, R., "Learning Optimal Control and Dynamical Structure of Global Trajectory Search Problems with Diffusion Models," Submitted, 10 2024. https://doi.org/10.48550/arXiv.org.2410.02976, URL https://arxiv.org/abs/2410.02976.

[54] Gräbner, J., and Beeson, R., "Semi-Supervised Global Search Capabilities for Optimal Low Thrust Spacecraft Trajectories over an Indirect Approach," *35th AAS/AIAA Space Flight Mechanics Meeting*, 2025.

[55] Li, A., Ding, Z., Dieng, A. B., and Beeson, R., "Efficient and Guaranteed-Safe Non-Convex Trajectory Optimization with Constrained Diffusion Model," *ICLR 2024 Workshop GenAI4DM*, 2024. https://doi.org/10.48550/arXiv.2403.05571, URL https://arxiv.org/abs/2403.05571v1.

[56] Jiang, Z., Zheng, Y., Tan, H., Tang, B., and Zhou, H., "Variational Deep Embedding: An Unsupervised and Generative Approach to Clustering," , 2017. https://doi.org/10.48550/arXiv.1611.05148, URL https://arxiv.org/abs/1611.05148v3.

[57] Xiong, L., Xu, K., Tian, K., Shao, Y., Tang, L., Gao, G., Zhang, M., Jiang, T., and Zhang, Q. C., "SCALE method for single-cell ATAC-seq analysis via latent feature extraction," *Nature Communications*, Vol. 10, No. 1, 2019, p. 4576. https://doi.org/10.1038/s41467-019-12630-7, URL https://doi.org/10.1038/s41467-019-12630-7.

[58] Sohn, K., Lee, H., and Yan, X., "Learning Structured Output Representation using Deep Conditional Generative Models," *Advances in Neural Information Processing Systems*, Vol. 28, edited by C. Cortes, N. Lawrence, D. Lee, M. Sugiyama, and R. Garnett, Curran Associates, Inc., 2015. URL https://proceedings.neurips.cc/paper_files/paper/2015/file/8d55a249e6baa5c06772297520da2051-Paper.pdf.

[59] Kingma, D. P., and Welling, M., "An Introduction to Variational Autoencoders," *Foundations and Trends® in Machine Learning*, Vol. 12, No. 4, 2019, pp. 307–392. https://doi.org/10.1561/2200000056, URL http://dx.doi.org/10.1561/2200000056.

[60] Hochreiter, S., and Schmidhuber, J., "Long Short-Term Memory," *Neural Computation*, Vol. 9, No. 8, 1997, pp. 1735–1780. https://doi.org/10.1162/neco.1997.9.8.1735, URL https://doi.org/10.1162/neco.1997.9.8.1735.

[61] Graves, A., and Schmidhuber, J., "Framewise phoneme classification with bidirectional LSTM and other neural network architectures," *Neural Networks*, Vol. 18, No. 5, 2005, pp. 602–610. https://doi.org/https://doi.org/10.1016/j.neunet.2005.06.042, URL https://www.sciencedirect.com/science/article/pii/S0893608005001206, iJCNN 2005.

[62] Paszke, A., Gross, S., Chintala, S., Chanan, G., Yang, E., DeVito, Z., Lin, Z., Desmaison, A., Antiga, L., and Lerer, A., "Automatic differentiation in PyTorch," *NIPS 2017 Workshop Autodiff*, Long Beach, California, 2017. URL https://openreview.net/forum?id=BJJsrmfCZ.

[63] Kingma, D. P., and Ba, J., "Adam: A Method for Stochastic Optimization," , 2017. https://doi.org/10.48550/arXiv.1412.6980, URL https://arxiv.org/abs/1412.6980v9.